\documentclass{amsart}
\usepackage{all,subfigure}

\newcommand{\AP}{\operatorname{AP}}
\newcommand{\VP}{\operatorname{VP}}

\newcommand{\AC}{\operatorname{Win}\nolimits}
\newcommand{\mF}{\mathcal{F}}
\newcommand{\Xu}{{X(\liea{u})}}
\newcommand{\tXu}{\tilde{X}(\liea{u})}

\newcommand{\Win}{\operatorname{Win}}

\begin{document}

\begin{abstract}
We completely describe the higher secant dimensions of all connected
homogeneous projective varieties of dimension at most $3$, in all
possible equivariant embeddings. In particular, we calculate these
dimensions for all Segre-Veronese embeddings of $\PP^1 \times \PP^1$,
$\PP^1 \times \PP^1 \times \PP^1$, and $\PP^2 \times \PP^1$, as well
as for the variety $\mF$ of incident point-line pairs in $\PP^2$. For
$\PP^2 \times \PP^1$ and $\mF$ the results are new, while the proofs for
the other two varieties are more compact than existing proofs. Our main
tool is the second author's tropical approach to secant dimensions.
\end{abstract}

\title{Secant dimensions of low-dimensional homogeneous varieties}
\author[K.~Baur]{Karin Baur}
\address[Karin Baur]{
Department of Mathematics\\
University of Leicester\\
University Road, Leicester LE1 7RH, U.K.}
\thanks{The first author is supported by EPSRC grant number
GR/S35387/01.}
\email{k.baur@mcs.le.ac.uk}

\author[J.~Draisma]{Jan Draisma}
\address[Jan Draisma]{
Department of Mathematics and Computer Science\\
Technische Universiteit Eindhoven\\
P.O. Box 513, 5600 MB Eindhoven, Netherlands}
\thanks{The second author is supported by DIAMANT, an NWO
mathematics cluster.}
\email{j.draisma@tue.nl}

\maketitle

\section{Introduction and results}

Let $K$ be an algebraically closed field of characteristic $0$; all
varieties appearing here will be over $K$. Let $G$ be a connected affine
algebraic group, and let $X$ be a projective variety on which $G$ acts
transitively. An {\em equivariant embedding} of $X$ is by definition
a $G$-equivariant injective morphism $\iota: X \rightarrow \PP(V)$,
where $V$ is a finite-dimensional (rational) $G$-module, subject to
the additional constraint that $\iota(X)$ spans $\PP(V)$.  The {$k$-th
(higher) secant variety} $k\iota(X)$ of $\iota(X)$ is the closure in
$\PP(V)$ of the union of all subspaces of $\PP(V)$ spanned by $k$ points
on $\iota(X)$. The {\em expected dimension} of $k\iota(X)$ is $\min\{k(\dim X +
1)-1,\dim V-1\}$; this is always an upper bound on $\dim k\iota(X)$. We
call $k\iota(X)$ {\em non-defective} if it has the expected dimension, and
{\em defective} otherwise. We call $\iota$ non-defective if $k\iota(X)$ is
non-defective for all $k$, and defective otherwise.

We want to compute the secant dimensions of $\iota(X)$ for all $X$ of
dimension at most $3$ and all $\iota$. This statement really concerns
only finitely many pairs $(G,X)$: Indeed, as $X$ is projective and
$G$-homogeneous, the stabiliser of any point in $X$ is parabolic (see
\cite[\S 11]{Borel91}) and therefore contains the solvable radical $R$
of $G$. But then $R$ also acts trivially on the span of $\iota(X)$,
which is $\PP(V)$, so that we may replace $G$ by the quotient $G/R$,
which is semisimple. In addition, we may and will assume that $G$ is
simply connected. Now $V$ is an irreducible $G$-module, and $\iota(X)$
is the unique closed orbit of $G$ in $\PP V$, the cone over which in $V$
is also known as the {\em cone of highest weight vectors}. Conversely,
recall that for two dominant weights $\lambda$ and $\lambda'$ the minimal
orbits in the corresponding projective spaces $\PP V(\lambda)$ and $\PP
V(\lambda')$ are isomorphic (as $G$-varieties) if and only if $\lambda$
and $\lambda'$ have the same support on the basis of fundamental
weights. So, to prove that all equivariant embeddings of a fixed $X$
are non-degenerate, we have to consider all possible dominant weights
with a fixed support.

Now there are precisely seven pairs $(G,X)$ with $\dim X \leq 3$, namely
$(\lieg{SL}_2^i, (\PP^1)^i)$ for $i=1,2,3$, $(\lieg{SL}_3,\PP^2)$,
$(\lieg{SL}_3 \times \lieg{SL}_2, \PP^2 \times \PP^1)$, $(\lieg{SL}_4,
\PP^3)$, and $(\lieg{SL}_3, \mF)$, where $\mF$ is the variety of flags
$p \subset l$ with $p,l$ a point and a line in $\PP^2$, respectively.
The equivariant embeddings of $\PP^i$ for $i=1,2,3$ are the Veronese
embeddings; their higher secant dimensions---and indeed, all higher
secant dimensions of Veronese embeddings of projective spaces of
arbitrary dimensions---are known from the work of Alexander and
Hirschowitz \cite{Alexander95}. In low dimensions there also exist
tropical proofs for these results: $\PP^1$ and $\PP^2$ were given as
examples in \cite{Draisma06a}, and for $\PP^3$ see the Master's thesis
of Silvia Brannetti \cite{Brannetti07}. The other varieties are covered
by the following theorems.

First, the equivariant embeddings of $\PP^1 \times \PP^1$ are
the Segre-Veronese embeddings, parametrised by the degree $(d,e)$
(corresponding to the highest weight $d \omega_1 +e \omega_2$ where the
$\omega_i$ are the fundamental weights), where we may assume $d \geq e$.
The following theorem is known in the literature; see for instance
\cite[Theorem 2.1]{Catalisano05c} and the references there.  Our proof
is rather short and transparent, and serves as a good introduction to
the more complicated proofs of the remaining theorems.

\begin{thm} \label{thm:P1P1}
The Segre-Veronese embedding of $\PP^1 \times \PP^1$ of degree $(d,e)$
with $d \geq e \geq 1$ is non-defective unless $e=2$ and $d$ is even, in which
case the $(d+1)$-st secant variety has codimension $1$ rather than the
expected $0$.
\end{thm}

The equivariant embeddings of $\PP^1 \times \PP^1 \times \PP^1$ and of
$\PP^2 \times \PP^1$ are also Segre-Veronese embeddings.  While writing
this paper we found out that the following theorem has already been
proved in \cite{Catalisano07}. We include our proof because we need its
building blocks for the other $3$-dimensional varieties.

\begin{thm} \label{thm:P1P1P1}
The Segre-Veronese embedding of $\PP^1 \times \PP^1 \times \PP^1$ of
degree $(d,e,f)$ with $d \geq e \geq f \geq 1$ is non-defective
unless 
\begin{enumerate}
\item $e=f=1$ and $d$ is even, in which case the $(d+1)$-st
secant variety has codimension $1$ rather than the expected
$0$, or 
\item $d=e=f=2$, in which case the $7$-th secant variety has
codimension $1$ instead of the expected $0$.
\end{enumerate}
\end{thm}

\begin{thm} \label{thm:P2P1}
The Segre-Veronese embedding of $\PP^2 \times \PP^1$ of
degree $(d,e)$ with $d,e \geq 1$ is non-defective unless
\begin{enumerate}
\item $d=2$ and $e=2k$ is even, in which case the $(3k+1)$-st secant
variety has codimension $3$ rather than the expected $2$ and
the $(3k+2)$-nd secant variety has codimension $1$ rather
than $0$; or 
\item $d=3$ and $e=1$, in which case the $5$-th secant
variety has codimension $1$ rather than the expected $0$.
\end{enumerate}
\end{thm}

Finally, the equivariant embeddings of $\mF$ are the minimal orbits
in $\PP V$ for any irreducible $\lieg{SL}_3$-representation of highest
weight $d \omega_1+e \omega_2$.

\begin{thm} \label{thm:F}
The image of $\mF$ in $\PP V$, for $V$ an irreducible
$\lieg{SL}_3$-representation of highest weight $d \omega_1+e\omega_2$
with $d,e \geq 1$ is non-defective unless
\begin{enumerate}
\item $d=e=1$, in which the $2$nd secant variety has codimension
$1$ rather than $0$, or 
\item $d=e=2$, in which the $7$th secant variety has
codimension $1$ rather than $0$. 
\end{enumerate}
\end{thm}

To the best of our knowledge Theorems \ref{thm:P2P1} and \ref{thm:F}
are new. Moreover, $\mF$ seems to be the first settled case where
maximal tori in $G$ do not have dense orbits. Our proofs of Theorems
\ref{thm:P1P1} and \ref{thm:P1P1P1} are more compact than their
original proofs \cite{Catalisano05c,Catalisano07}. Moreover, the planar
proof of Theorem \ref{thm:P1P1} serves as a good introduction to the
more complicated induction in the three-dimensional cases, while parts
of the proof of Theorem \ref{thm:P1P1P1} are used as building blocks in
the remaining proofs.

We will prove our theorems using a polyhedral-combinatorial lower
bound on higher secant dimensions introduced by the second author in
\cite{Draisma06a}. Roughly this goes as follows: to a given $X$ and
$V$ we associate a finite set $B$ of points in $\RR^{\dim X}$, which
parametrises a certain basis in $V$. Now to find a lower bound on $\dim
kX$ we {\em maximise}
\[ \sum_{i=1}^k [1+\dim \Aff_\RR \Win_i(f)] \]
over all $k$-tuples $f=(f_1,\ldots,f_n)$ of affine-linear functions on
$\RR^{\dim X}$, where $\Win_i(f)$ is the set of points in $B$ where $f_i$
is strictly smaller than all $f_j$ with $j \neq i$, and where $\Aff_\RR$
denotes taking the affine span. Typically, this maximum equals $1$ plus
the expected dimension of $\dim k\iota(X)$, and then we are done. If
not, then we need other methods to prove that $k\iota(X)$ is indeed
defective---but most defective cases above are known in the literature.

\begin{figure}
\includegraphics[scale=.4]{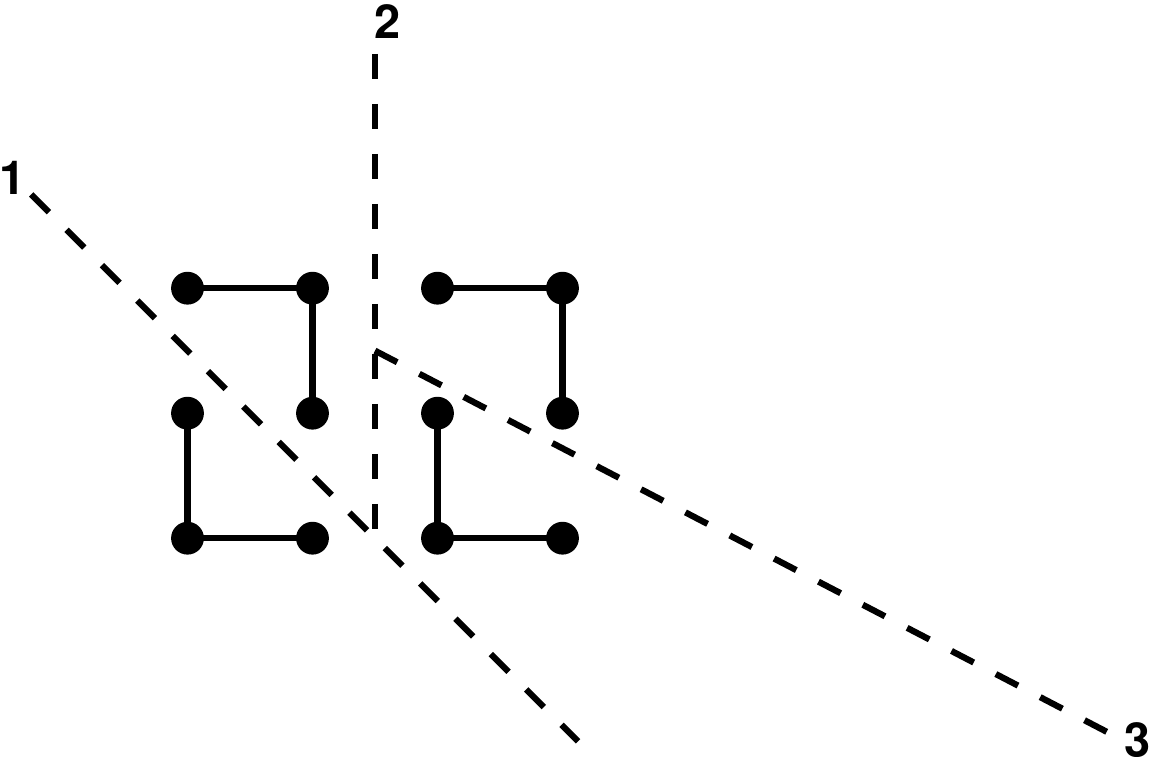}
\caption{The embedding $\Seg \circ (v_3 \times v_2)$ of
$\PP^1 \times \PP^1$ is non-defective.}
\label{fig:sl2sl2_32a}
\end{figure}

This optimisation problem may sound somewhat far-fetched, so as a
motivation we now carry out our proof in one particular case: For the
Segre-Veronese embedding of $X=\PP^1 \times \PP^1$ of degree $(d,e)$
the set $B$ is the grid $\{0,\ldots,d\} \times \{0,\ldots,e\} \subseteq
\RR^2$. Take for instance $d=3$ and $e=2$. In Figure \ref{fig:sl2sl2_32a}
the points in $B$ are grouped into four triples spanning the plane. It
is easy to see---see Lemma \ref{lm:Relation} below---that there exist
affine-linear $f_1,\ldots,f_4$ inducing this partition, so that $4X$
has the expected dimension $4 \cdot 3 - 1=11$.

Our tropical approach is conceptually very simple, and closely
related to Sturmfels-Sullivant's combinatorial secant varieties
\cite{Sturmfels05}, Miranda-Dumitrescu's degeneration approach (private
communication), and Develin's tropical secant varieties of linear spaces
\cite{Develin06}. What we find very surprising and promising is that
strong results on secant varieties of non-toric varieties such as $\mF$
can be proved with our approach. 

The remainder of this paper is organised as follows. In Section
\ref{sec:Bounds} we recall the tropical approach, and prove a lemma
that will help us deal with the flag variety. The tropical approach
depends rather heavily on a parameterisation of $X$, and in Section
\ref{sec:PolynomialMap} we introduce the polynomial maps that we will
use. In particular, we give, for any minimal orbit (not necessarily of
low dimension, and not necessarily toric), a polynomial paramaterisation
whose tropicalisation has an image of the right dimension; these
tropical parameterisations are useful in studying tropicalisations
of minimal orbits; see Remark \ref{re:TropicalMinOrbit}. Finally,
Sections \ref{sec:P1P1}--\ref{sec:F} contain the proofs of Theorems
\ref{thm:P1P1}--\ref{thm:F}, respectively.

\section{The tropical approach} \label{sec:Bounds}

\subsection{Two optimisation problems}
\label{ssec:Optimisation}

We recall from \cite{Draisma06a} a polyhedral-combinatorial optimisation
problem that plays a crucial role in the proofs of our theorems; here
$\AP$ abbreviates {\em Affine Partition}.

\begin{prb}[$\AP(A,k)$]
Let $A=(A_1,\ldots,A_n)$ be a sequence of finite subsets of $\RR^m$ and
let $k \in \NN$. For any $k$-tuple $f=(f_1,\ldots,f_k)$ of affine-linear
functions on $\RR^m$ let the sets $\AC_i(f),\ i=1,\ldots,k$ be defined
as follows. For $b=1,\ldots,n$ we say that $f_i$ {\em wins} $b$ if $f_i$
attains its minimum on $A_b$ in a unique $\alpha \in A_b$, and if this
minimum is strictly smaller than all values of all $f_j, j \neq i$
on $A_b$. The vector $\alpha$ is then called a {\em winning direction}
of $f_i$. Let $\AC_i(f)$ denote the set of winning directions of $f_i$.

\text{\bf Maximise } $\sum_{i=1}^k [1+\dim \Aff_\RR \AC_i(f)]$ over all
$k$-tuples $f$ of affine-linear functions on $\RR^m$; call the maximum
$\AP^*(A,k)$.
\end{prb}

Note that if every $A_b$ is a singleton $\{\alpha_b\}$, then
$\AC_i(f)$ is just the set of all $\alpha_b$ on which $f_i$ is
smaller than all other $f_j, j \neq i$. We will then also write
$\AP(\{\alpha_1,\ldots,\alpha_n\},k)$ for the optimisation problem
above. In this case we are really optimising over all possible {\em
regular subdivisions} of $\RR^m$ into $k$ open cells.  Each such
subdivision induces a partition of the $\alpha_b$ into the sets $\AC_i(f)$
(at least if no $\alpha_b$ lies on a border between two cells, but this
is easy to achieve without decreasing the objective function). As it
is sometimes hard to imagine the existence of affine-linear functions
inducing a certain regular subdivision of space, we have the following
observation, due to Immanuel Halupczok and the second author.

\begin{lm} \label{lm:Cutting}
Let $S$ be a finite set in $\RR^m$, let $f_1,\ldots,f_k$ be affine-linear
functions on $\RR^m$, and let $g_1,\ldots,g_l$ also be affine-linear
functions on $\RR^m$. Let $S_i$ be the subset of $S$ where $f_i<f_j$ for
all $j \neq i$, and let $T_i$ be the subset of $S_1$ where $g_i<g_j$ for
all $j \neq i$. Then there exist affine-linear functions $h_1,\ldots,h_l$
such that 
\begin{enumerate} 
\item $h_i < h_j$ on $T_i$ for $i,j=1,\ldots,l$ and $j \neq i$;
\item $h_i < f_j$ on $T_i$ for $i=1,\ldots,l$ and $j=2,\ldots,k$; and 
\item $f_i < h_j$ on $S_i$ for $i=2,\ldots,k$ and
$j=1,\ldots,l$.
\end{enumerate}
In other words, the functions $h_1,\ldots,h_l,f_2,\ldots,f_k$ together
induce the partition $T_1,\ldots,T_l,S_2,\ldots,S_k$ of $S$.
\end{lm}

\begin{proof}
Take $h_i=f_1 + \epsilon g_i$ for $\epsilon$ positive and sufficiently
small.
\end{proof}

This lemma implies, for instance, that one may find appropriate
$\AC_i(f)$ (still for the case of singletons $A_b$) by repeatedly
cutting polyhedral pieces of space in half. For instance, in Figure
\ref{fig:sl2sl2_32a} the plane is cut into four pieces by three straight
cuts. Although this is not a regular subdivision of the plane, by the
lemma there does exist a regular subdivision of the plane inducing the
same partition on the $12$ points. 

Lemma \ref{lm:Cutting} can only immediately be applied to $\AP$ if
the $A_b$ are singletons, while the $A_b$ in our application to the
$3$-dimensional flag variety $\mF$ are not. We get around this difficulty
by giving a lower bound on $\AP^*(A,k)$ for more general $A$ in terms of
$\AP^*(A',k)$ for some sequence $A'$ of singletons.  In the following
lemmas a {\em convex polyhedral cone} in $\RR^m$ is by definition the
set of nonnegative linear combinations of a finite set in $\RR^m$, and
it is called {\em strictly convex} if it does not contain any non-trivial
linear subspace of $\RR^m$.

\begin{lm} \label{lm:Cone}
Let $A=(\{\alpha_1\},\ldots,\{\alpha_n\})$ be an $n$-tuple of singleton
subsets of $\RR^m$. Furthermore, let $k \in \NN$, let $Z$ be a strictly
convex polyhedral cone in $\RR^m$, and let $f$ be a $k$-tuple of
affine-linear functions on $\RR^m$. Then the value of $\AP(A,k)$ at $f$
is also attained at some $f'=(f_1',\ldots,f_k')$ for which every $f_i'$
is strictly decreasing in the $z$-direction, for every $z \in Z\setminus
\{0\}$.
\end{lm}

\begin{proof}
By the strict convexity of $Z$, there exists a linear function $f_0$
on $\RR^m$ such that every $f_j+f_0$ is strictly decreasing in the
$z$-direction, for every $z \in Z$. But since
\[ f_i(\alpha)<f_j(\alpha) \Leftrightarrow f_i(\alpha) +
f_0(\alpha)<f_j(\alpha) + f_0(\alpha) \]
we have $\AC_i((f_j+f_0)_j)=\AC_i(f)$ for all $i$, and we
are done.
\end{proof}

It is crucial in this proof that only values of $f_i$ and $f_j$ at {\em the
same} $\alpha$ are compared---that is why we have restricted ourselves
to singleton-$\AP$ here.

\begin{lm} \label{lm:Relation}
Let $A=(A_1,\ldots,A_n)$ be a $k$-tuple of finite subsets of $\RR^m$
and let $k \in \NN$. Furthermore, let $Z$ be a strictly convex polyhedral
cone in $\RR^m$ and define a partial order $\leq$ on $\RR^m$ by
\[ p \leq q :\Leftrightarrow p-q \in Z. \]
Suppose that for every $b$, $A_b$ has a unique minimal element $\alpha_b$
with respect to this order. Then we have
\[ \AP^*(A,k) \geq \AP^*(\{\alpha_1,\ldots,\alpha_n\},k) \]
\end{lm}

\begin{proof}
Let $d^*=\AP^*(\{\alpha_1,\ldots,\alpha_n\},k)$. By Lemma \ref{lm:Cone} there
exists a $k$-tuple $f=(f_1,\ldots,f_k)$ of affine-linear functions on
$\RR^m$ for which $\AP(\{\alpha_1,\ldots,\alpha_n\},k)$ also has value $d^*$ and for
which every $f_i$ is strictly decreasing in all directions in $Z$. We
claim that the value of $\AP(A,k)$ at this $f$ is also $d^*$. Indeed,
fix $b \in B$ and consider all $f_i(\alpha)$ with $\alpha \in A_b$ and
$i=1,\ldots,k$. Because $\alpha_b-\alpha \in Z$ for all $\alpha \in A_b$
and because every $f_i$ is strictly decreasing in the directions in $Z$,
we have $f_i(\alpha_b) < f_i(\alpha)$ for all $\alpha \in A_b \setminus
\{\alpha_b\}$ and all $i$. Hence the minimum, over all pairs $(i,\alpha)
\in \{1,\ldots,k\} \times A_b$, of $f_i(\alpha)$ can only be attained
in pairs for which $\alpha=\alpha_b$.  Therefore, in computing the
value at $f$ of $\AP(A,k)$ the elements of $A_b$ unequal to $\alpha_b$
can be ignored. We conclude that $\AP(A,k)$ has value $d^*$ at $f$,
as claimed. This shows the inequality.
\end{proof}

\subsection{Tropical bounds on secant dimensions}
\label{ssec:TropicalBounds}

Rather than working with projective varieties, we work with the
affine cones over them. So suppose that $C \subseteq K^n$ is a closed
cone (i.e., closed under scalar multiplication with $K$), and set 
\[ kC:=\overline{\{v_1+\ldots+v_k \mid v_1,\ldots,v_k \in C\}}. \] 
Suppose that $C$ is unirational, and choose a polynomial map
$f=(f_1,\ldots,f_n):K^m \rightarrow C \subseteq K^n$ that maps $K^m$
dominantly into $C$.  Let $x=(x_i)_{i=1}^m$ and $y=(y_b)_{b=1}^n$ be the
standard coordinates on $K^m$ and $K^n$. The tropical approach depends
very much on coordinates; in particular, one would like $f$ to be sparse.
For every $b=1,\ldots,n$ let $A_b$ be the set of $\alpha \in \NN^m$
for which the monomial $x^\alpha$ has a non-zero coefficient in $f_b$,
and set $A:=(A_1,\ldots,A_n)$.

\begin{thm}[\cite{Draisma06a}] \label{thm:TropicalBounds}
For all $k \in \NN$, $\dim kC \geq \AP^*(A,k)$.
\end{thm}

\begin{re}
In fact, in \cite{Draisma06a} this is proved provided that $\bigcup_b A_b$ is
contained in an affine hyperplane not through $0$, but this can always
be achieved by taking a new map $f'(t,x):=t f(x)$ into $C$, without
changing the optimisation problem $\AP(A,k)$.
\end{re}

In Section \ref{sec:PolynomialMap} we introduce a polynomial map $f$ for
general minimal orbits that seems suitable for the tropical approach, and
after that we specialise to low-dimensional varieties under
consideration. 

\subsection{Non-defective pictures}
Our proofs will be entirely pictorial: given a set $B$ of lattice points
in $\ZZ^2$ or $\ZZ^3$ according as $\dim X=2$ or $\dim X=3$, we solve the
optimisation problem $\AP(B,k)$ for all $k$. To this end, we will exhibit
a partition of $B$ into parts $B_i$ such that there exist affine-linear
functions $f_i$ on $\RR^2$ or $\RR^3$, exactly one for each part, with
$B_i=\AC_i(f)$. If each $B_i$ is affinely independent, and if moreover
the affine span of each $B_i$ has $\dim X$, except possibly for one
single $B_i$, then we call the picture {\em non-defective}, as it shows,
by Theorem \ref{thm:TropicalBounds}, that all secant varieties of $X$
in the given embedding have the expected dimension.

\begin{figure}
\includegraphics[scale=.5]{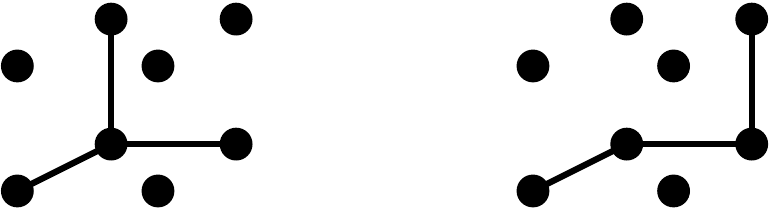}
\caption{Types $1$ (left) and $2$ (right).}
\label{fig:configs}
\end{figure}

The full-dimensional $B_i$ that we will use will have very simple
shapes: in dimension $2$ they will all be equivalent, up to integral
translations and rotations over multiples of $\pi/2$, to the triple
$\{0,e_1,e_2\}$. In dimension $3$ they will almost all be equivalent
to either $\{0,e_1,e_2,e_3\}$ (type $1$) or $\{e_1,0,e_2,e_2+e_3\}$
(type $2$), up to Euclidean transformations preserving the lattice;
see Figure \ref{fig:configs}. Only in case of the flag-variety $\mF$
we will occasionally use more general pictures. 

The $f_i$ will not be explicitly computed. Indeed, in all cases
their existence follows from a tedious but easy application of Lemma
\ref{lm:Cutting}: one can repeatedly cut $\RR^2$ or $\RR^3$ into pieces
by affine hyperplanes, such that eventually the desired partition of $B$
into the $B_i$ is attained. For instance, in Figure \ref{fig:sl2sl2_32a}
three cuts, labelled $1,2,3$ consecutively, give the desired partition
of the twelve points. 

\section{A polynomial map} \label{sec:PolynomialMap}
We retain the setting of the Introduction: $G$ is a simply connected,
connected, semisimple algebraic group, $V$ is a $G$-module, and we wish
to determine the secant dimensions of $X$, the unique closed orbit of $G$
in $\PP V$. Let $C$ be the affine cone in $V$ over $X$.
Fix a Borel subgroup $B$ of $G$, let $T$ be a maximal torus of $B$
and let $v_\lambda \in V$ span the unique $B$-stable one-dimensional
subspace of $V$; $\lambda$ denotes the $T$-weight of $v_\lambda$, i.e.,
the highest weight of $V$. Let $P \supseteq B$ be the stabiliser in $G$
of $Kv_\lambda$ (so that $X \cong G/P$ as a $G$-variety) and
let $U$ be the unipotent radical of the parabolic subgroup opposite to
$P$ and containing $T$. Let $\liea{u}$ denote the Lie algebra of $U_-$,
let $\Xu$ be the set of $T$-roots on $\liea{u}$, and set $\tXu:=\Xu \cup
\{0\}$. For every $\beta \in \Xu$ choose a vector $X_\beta$ spanning the
root space $\liea{u}_\beta$.  Moreover, fix an order on $\Xu$. Then it
is well-known that the polynomial map
\[ \Psi:K^{\tXu} \rightarrow V, 
	t \mapsto t_0 \prod_{\beta \in \Xu} \exp(t_\beta X_\beta)
	v_\lambda, \]
where the product is taken in the fixed order, maps dominantly into $C$.
This map will play the role of $f$ from Subsection
\ref{ssec:TropicalBounds}.

In what follows we will need the following notation: Let $X_\RR:=\RR
\otimes_\ZZ X(T)$ be the real vector space spanned by the character group
of $T$, let $\xi:\RR^{\Xu} \mapsto X_\RR$ send $r$ to $\sum_\beta r_\beta
\beta$ and also use $\xi$ for the map $\RR^{\tXu} \rightarrow X_\RR$
with the same definition; in both cases we call $\xi(r)$ the {\em weight}
of $r$.

Now for a basis of $V$: by the PBW-theorem, $V$ is the linear span of
all elements of the form $m_r:=\prod_{\beta \in \Xu} X_\beta^{r_\beta}
v_\lambda$ with $r \in \NN^{\Xu}$; the product is taken in the same
fixed order as before. Slightly inaccurately, we will call the $m_r$
PBW-monomials. Note that the $T$-weight of $m_r$ is $\lambda+\xi(r)$.  Let
$M$ be the subset of all $r \in \NN^{\Xu}$ for which $m_r$ is non-zero;
$M$ is finite. Let $B$ be a subset of $M$ such that $\{m_r \mid r \in B\}$
is a basis of $V$; later on we will add further restrictions on $B$. For
$b \in B$ let $\Psi_b$ be the component of $\Psi$ corresponding to $b$;
it equals $t_0$ times a polynomial in the $t_\beta,\ \beta \in \Xu$.
Let $A_b \subseteq \NN^{\Xu}$ denote the set of exponent vectors of
monomials having a non-zero coefficient in $\Psi_b / t_0$.

\begin{lm} \label{lm:Mons}
For $b_0 \in B$
\begin{enumerate}
\item $A_{b_0} \subseteq \{r \in M \mid \xi(r)=\xi(b_0)\}$,
and
\item $A_{b_0} \cap B=\{b_0\}$.
\end{enumerate}
\end{lm}

\begin{proof}
Expand $\Psi(t)/t_0$ as a linear combination of PBW-monomials:
\[ \Psi(t)/t_0=\sum_{r \in \NN^{\Xu}}
        \frac{t^r}{\prod_{\beta \in \Xu}(r_\beta!)} m_r. \]
So $t^r$ appears in $\Psi_{b_0}/t_0$ if and only if $m_r$ has a
non-zero $m_{b_0}$-coefficient relative to the basis $(m_b)_{b \in
B}$. Hence the first statement follows from the fact that every $m_r$
is a linear combination of the $m_b$ of the same $T$-weight as $m_r$,
and the second statement reflects the fact that for all $b_1 \in B$,
$m_{b_1}$ has precisely one non-zero coefficient relative to the basis
$(m_b)_{b \in B}$, namely that of $m_{b_1}$.
\end{proof}

Now Theorem \ref{thm:TropicalBounds} implies the following proposition.

\begin{prop} \label{prop:BoundAP}
$\dim kC \geq \AP^*((A_b)_{b \in B},k)$
\end{prop}

For Segre products of Veronese embeddings every $A_b$ is a singleton,
and we can use our hyperplane-cutting procedure immediately. For the
flag variety $\mF$ we will use Lemma \ref{lm:Relation} to bound $\AP^*$
by a singleton-$\AP^*$.

\begin{re} \label{re:TropicalMinOrbit}
To see that Proposition \ref{prop:BoundAP} has a chance of being useful,
it is instructive to verify that $\AP^*((A_b)_{b \in B},1)$ is, indeed,
$\dim C$, at least for some choices of $B$. Indeed, recall that $|\Xu|+1$
vectors $v_\lambda$ and $X_\beta v_\lambda, \beta \in \Xu$ are linearly
independent, so that we can take $B$ to contain the corresponding
exponent vectors, that is, $0$ and the standard basis vectors $e_\beta$
in $\NN^\Xu$. Now let $f_1:\RR^{\Xu} \rightarrow \RR$ send $r$ to
$\sum_{\beta \in \Xu} r_\beta$.  We claim that $\AP((A_b)_{b \in B},k)$
has value $\dim C$ at $(f_1)$. Indeed, $A_0=\{0\}$ and for every $b \in B$
of the form $e_\beta, \beta \in \Xu$ the set $A_b$ consists of $e_\beta$
itself, with $f_1$-value $1$, and exponent vectors having a $f_1$-value a
natural number $>1$. Hence $\AC_1(f_1)$ contains all $e_\beta$ and $0$---and
therefore spans an affine space of dimension $\dim C-1=\dim X$.

This observation is of some independent interest for tropical geometry:
going through the theory in \cite{Draisma06a}, it shows that the image
of the tropicalisation of $\Psi$ in the tropicalisation of $C$ has the
right dimension; this is useful in minimal orbits such as Grassmannians.
\end{re}

\section{Secant dimensions of $\PP^1 \times \PP^1$}
\label{sec:P1P1}

We retain the notation of Section \ref{sec:PolynomialMap}. To prove
Theorem \ref{thm:P1P1}, let $X=\PP^1 \times \PP^1$, $G=\lieg{SL}_2 \times
\lieg{SL}_2$, and $V=S^d(K^2) \otimes S^e(K^2)$. The polynomial map
\[ \Psi:(t_0,t_1,t_2) \mapsto t_0 (x_1+t_1 x_2)^d \otimes (x_1+t_2 x_2)^e, \] 
is dominant into the cone $C$ over $X$, and $M=B$ is the rectangular grid
$\{0,\ldots,d\} \times \{0,\ldots,e\}$. We may assume that $d \geq e$.

First, if $e=2$ and $d$ is even, then $(d+1)C$ is known to be defective,
that is, it does not fill $V$ but is given by some determinantal equation;
see \cite[Example 3.2]{Catalisano06}. The argument below will show that
its defect is not more than $1$.

\begin{figure}
\subfigure[$(d,e)=(1,1)$]
	{\includegraphics[scale=.4]{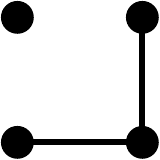}\label{fig:sl2sl2_11}}
	\hfill
\subfigure[$(d,e)=(2,1)$]
	{\includegraphics[scale=.4]{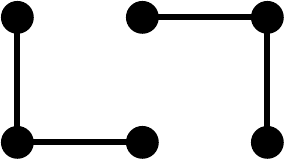}\label{fig:sl2sl2_21}}
	\hfill
\subfigure[$(d,e)=(3,1)$]
	{\includegraphics[scale=.4]{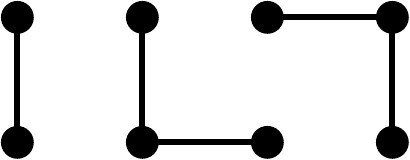}\label{fig:sl2sl2_31}}
	\hfill
\subfigure[$e=1$; induction]
	{\includegraphics[scale=.4]{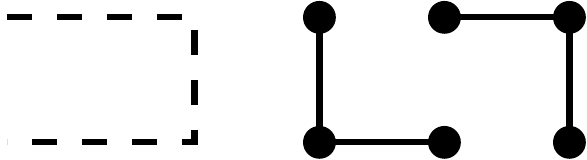}\label{fig:sl2sl2_1ind}}

	\hfill
\subfigure[$(d,e)=(2,2)$]
	{\includegraphics[scale=.4]{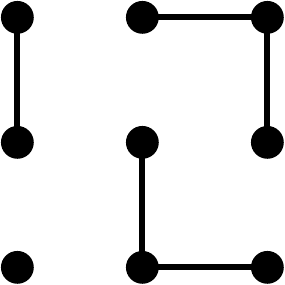}\label{fig:sl2sl2_22}}
	\hfill
\subfigure[$(d,e)=(3,2)$]
	{\includegraphics[scale=.4]{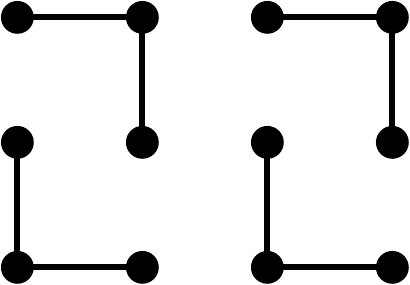}\label{fig:sl2sl2_32}}
	\hfill
\subfigure[$e=2$; induction]
	{\includegraphics[scale=.4]{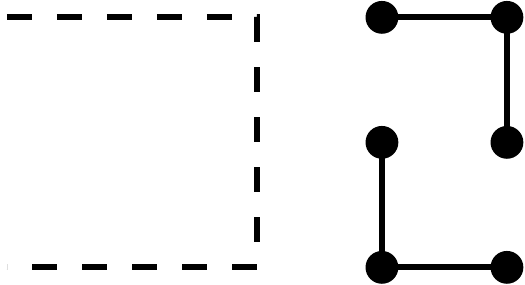}\label{fig:sl2sl2_2ind}}

	\hfill
\subfigure[$(d,e)=(3,3)$]
	{\includegraphics[scale=.4]{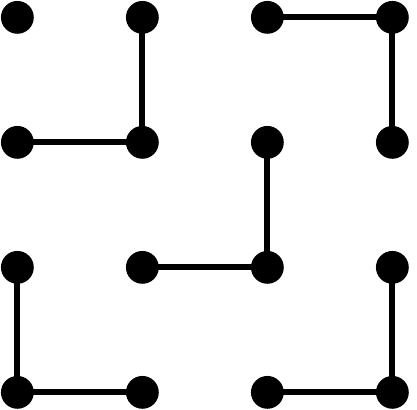}\label{fig:sl2sl2_33}}
	\hfill
\subfigure[$(d,e)=(4,3)$]
	{\includegraphics[scale=.4]{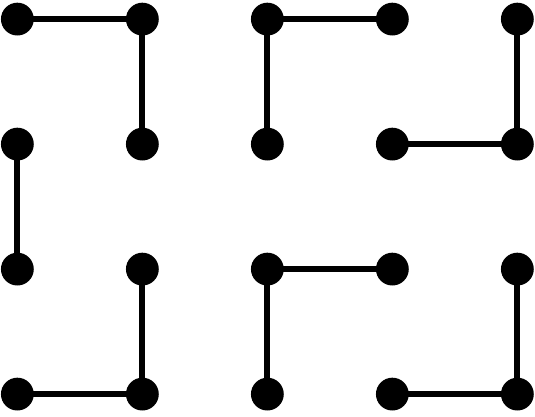}\label{fig:sl2sl2_43}}
	\hfill
\subfigure[$(d,e)=(5,3)$]{\includegraphics[scale=.4]{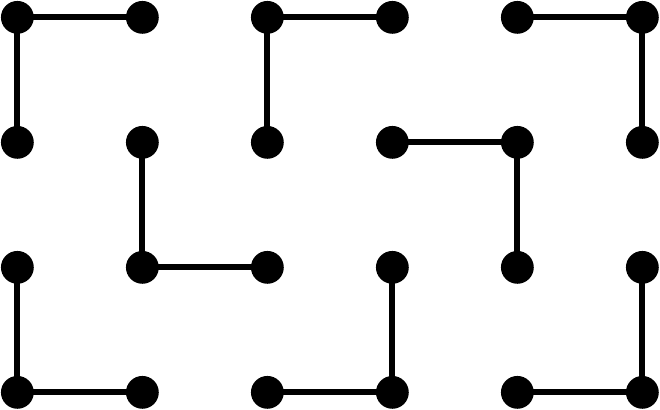}\label{fig:sl2sl2_53}}
	\hfill
\subfigure[$e=3$; induction]
	{\includegraphics[scale=.4]{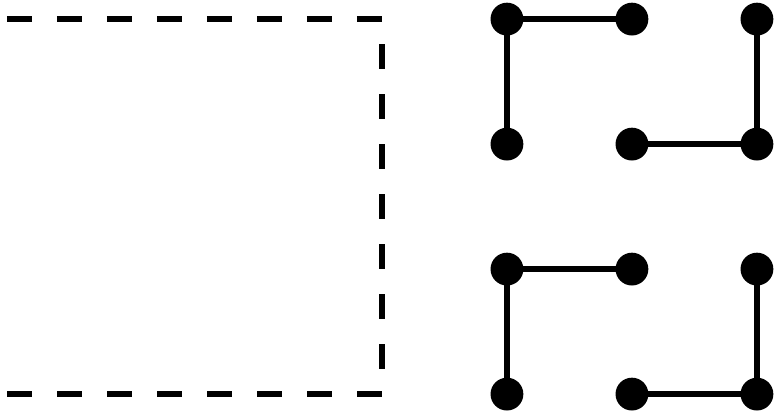}\label{fig:sl2sl2_3ind}}

	\hfill
\subfigure[$(d,e)=(3,4)$]
	{\includegraphics[scale=.4]{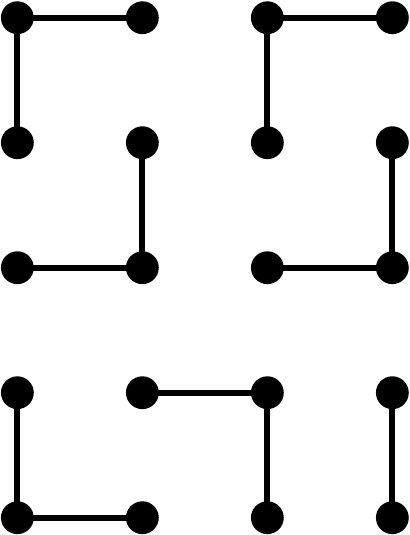}\label{fig:sl2sl2_34}}
	\hfill
\subfigure[$(d,e)=(4,4)$]
	{\includegraphics[scale=.4]{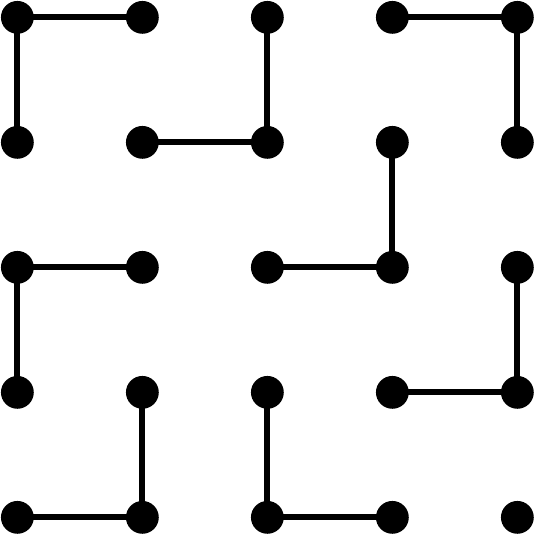}\label{fig:sl2sl2_44}}
	\hfill
\subfigure[$e=4$, induction]
	{\includegraphics[scale=.4]{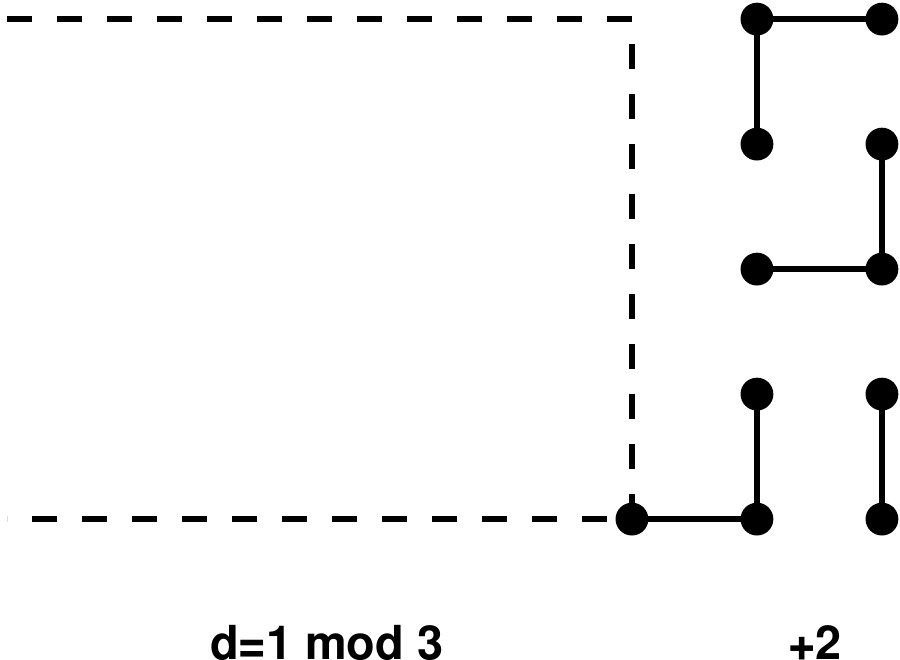}
	\includegraphics[scale=.4]{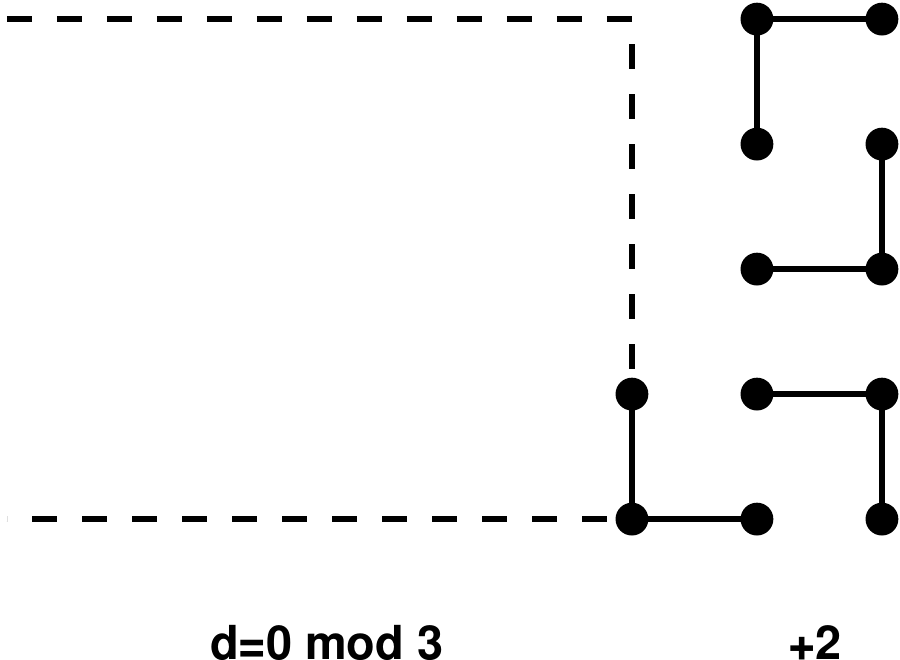}
	\includegraphics[scale=.4]{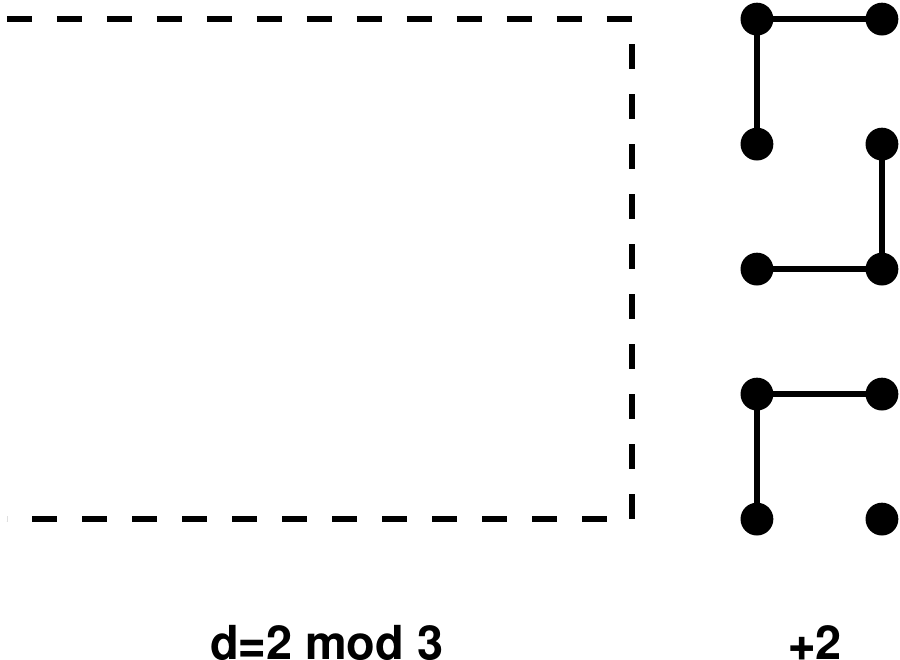}
	\label{fig:sl2sl2_4ind}
	}
	\hfill

\subfigure[$e=5$; induction]
	{\includegraphics[scale=.4]{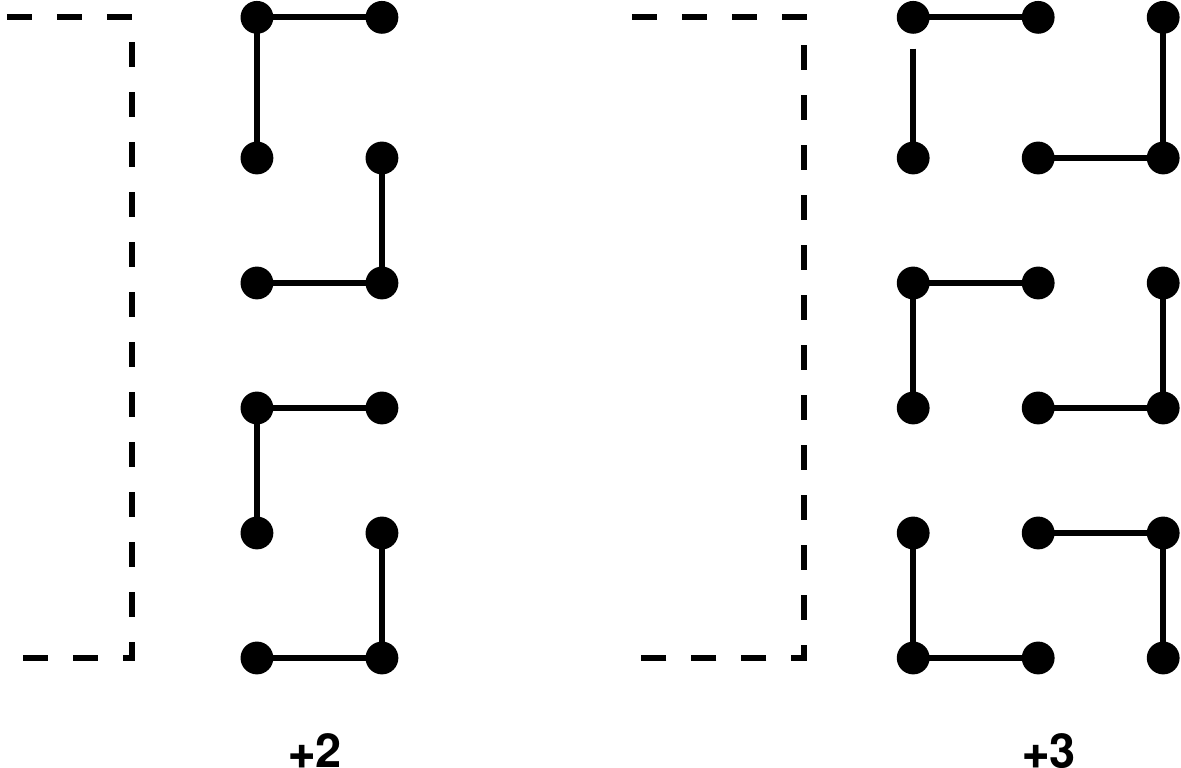}\label{fig:sl2sl2_5ind}}
\end{figure}

\begin{figure}
\subfigure[$(d,e)=(6,6)$]
	{\includegraphics[scale=.4]{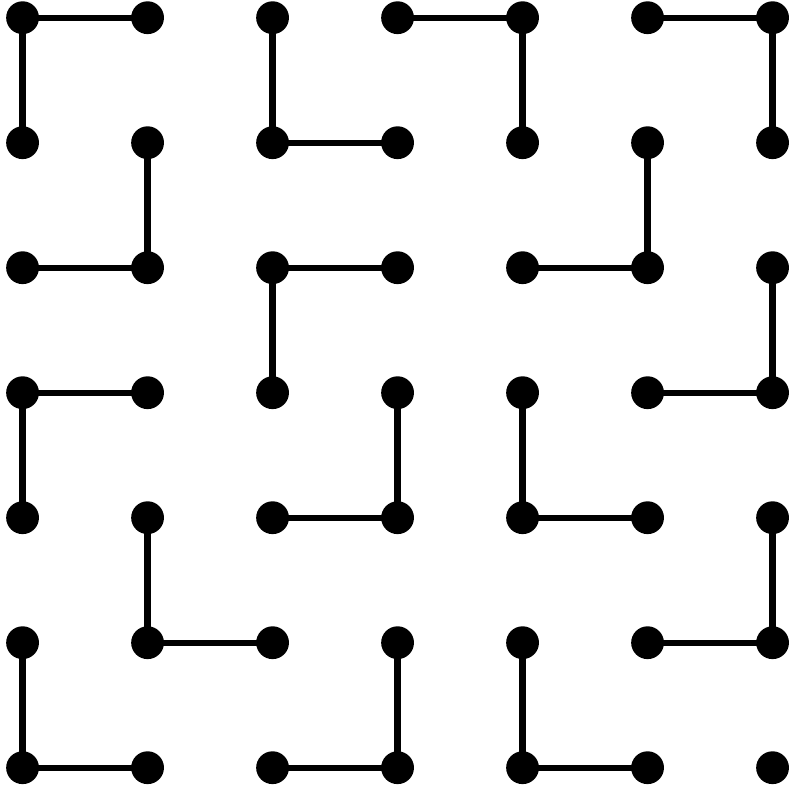}\label{fig:sl2sl2_66}}
	\hfill
\subfigure[$e=6$; induction]
	{\includegraphics[scale=.4]{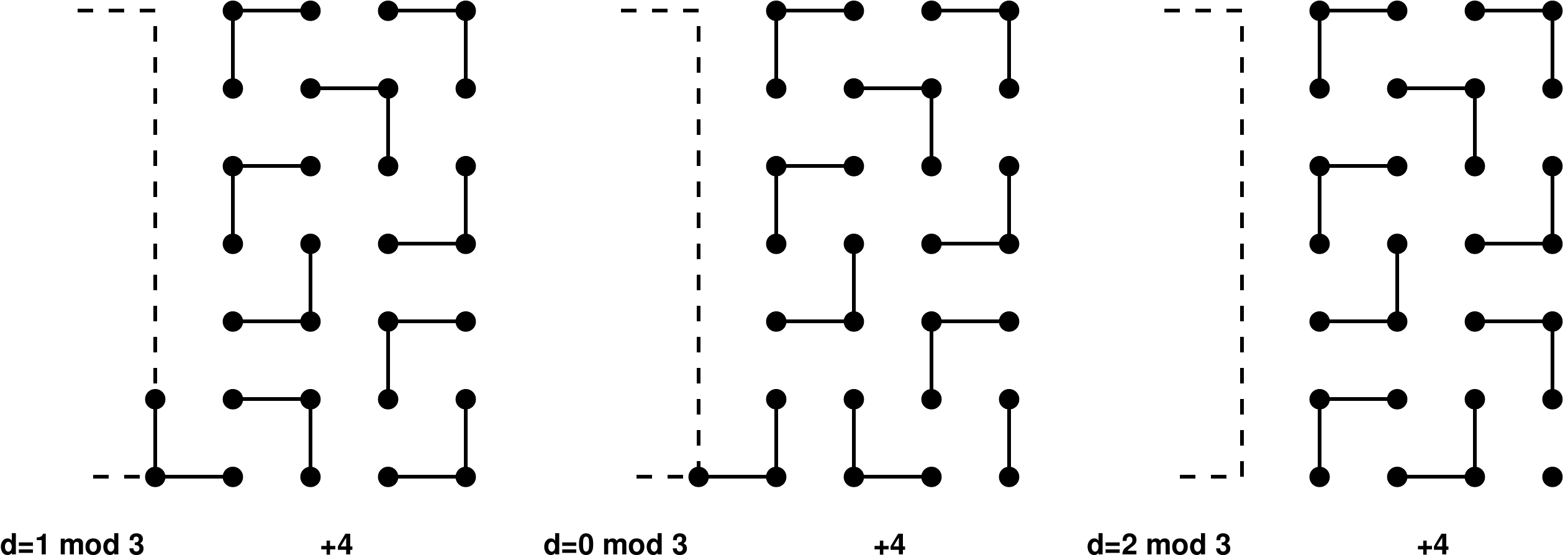}\label{fig:sl2sl2_6ind}}
	\hfill
\subfigure[$e=8$; induction]
	{\includegraphics[scale=.4]{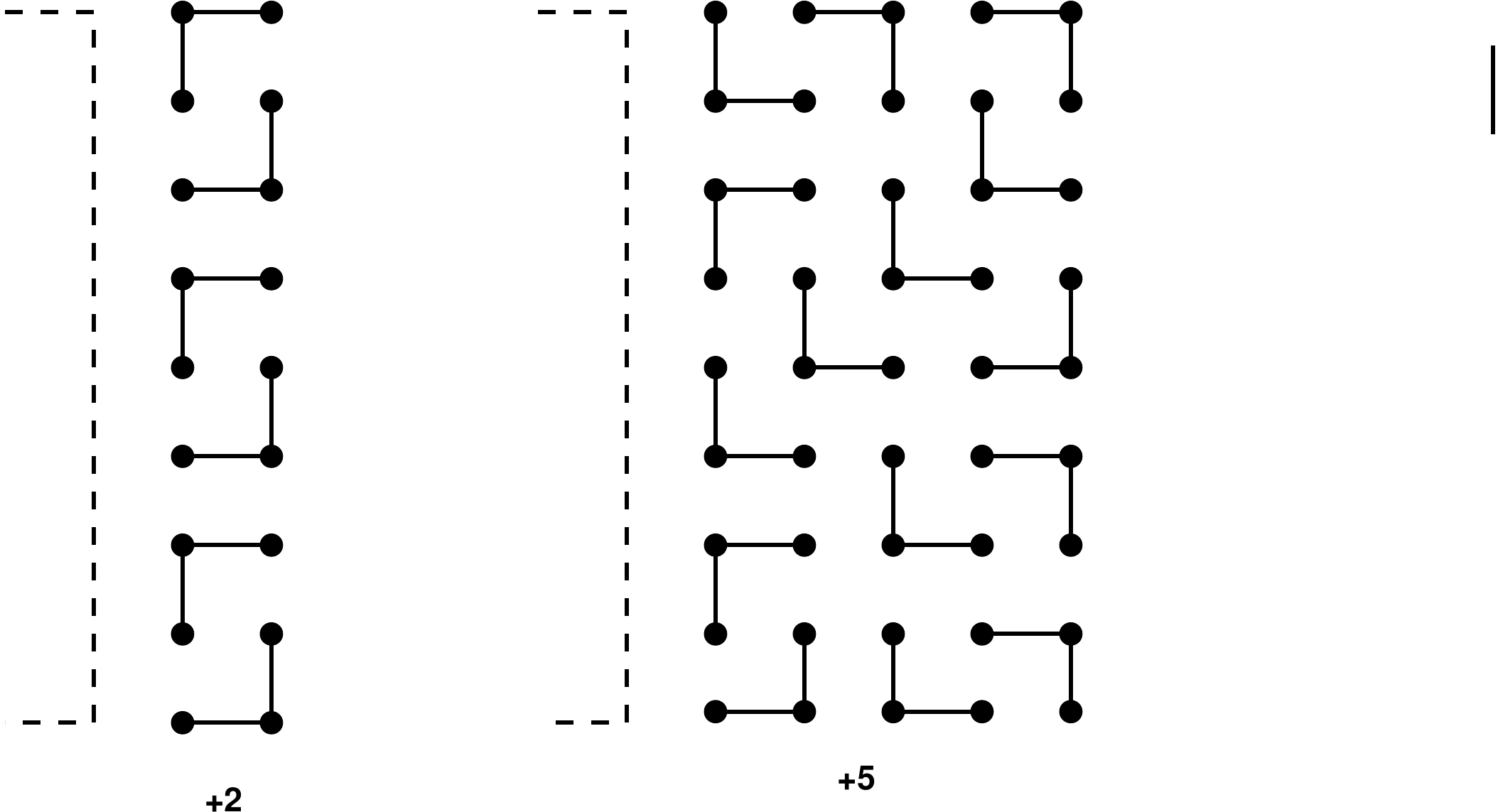}\label{fig:sl2sl2_8ind}}
\caption{More non-defective pictures for $\PP^1 \times \PP^1$.}
\label{fig:P1P1}
\end{figure}

Figure \ref{fig:P1P1} gives non-defective pictures for $e=1,2,3,4,5$
and $d \geq e$, except for $e=2$ and $d$ even. This implies, by
transposing pictures, that there exist non-defective pictures for
$e=6$ and $d=1,3,4,5$. Figure \ref{fig:sl2sl2_66} gives a non-defective
picture for $(d,e)=(6,6)$. Then, using the two induction steps in Figure
\ref{fig:sl2sl2_6ind}, we find non-defective pictures for $e=6$ and
all $d \neq 2$. A similar reasoning gives non-defective pictures for
$e=8$ and all $d \neq 2$. Finally, let $d \geq e \geq 6$ be arbitrary
with $(d,e) \not \in 2\NN \times \{2\}$. Write $e+1=6q+r$ with $r \in
\{0,2,4,5,7,9\}$. Then we find a non-defective picture for $(d,e)$ by
gluing $q$ non-defective pictures for $(d,5)$ and, if $r \neq 0$, one
non-defective picture for $(d,r-1)$ on top of each other. This proves
Theorem \ref{thm:P1P1}.

\clearpage

\section{Secant dimensions of $\PP^1 \times \PP^1 \times
\PP^1$} \label{sec:P1P1P1}

Now we turn to Theorem \ref{thm:P1P1P1}. Cutting to the chase,
$M=B$ is the block $\{0,\ldots,d\} \times \{0,\ldots,e\} \times
\{0,\ldots,f\}$. When convenient, we assume that $d \geq e \geq f$.
First, for $e=f=1$ and $d$ even, the $d+1$-st secant variety, which
one would expect to fill the space, is in fact known to be defective,
see \cite{Catalisano07}. The pictures below will show that the defect
is not more than $1$.

\begin{figure}
\subfigure[$(1,1,1)$]{\includegraphics[scale=.4]{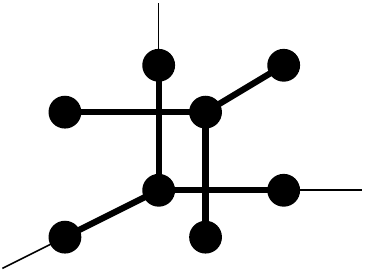}\label{fig:sl2sl2sl2_111}}
\hfill
\subfigure[$(2,1,1)$]{\includegraphics[scale=.4]{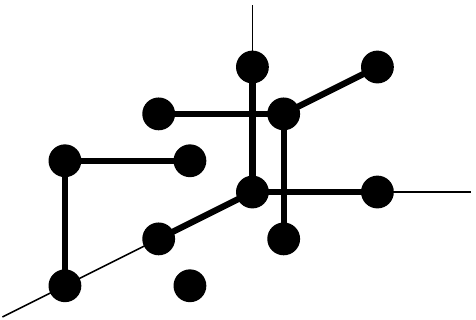}\label{fig:sl2sl2sl2_211}}
\hfill
\subfigure[induction for $(*,1,1)$]
	{\includegraphics[scale=.4]{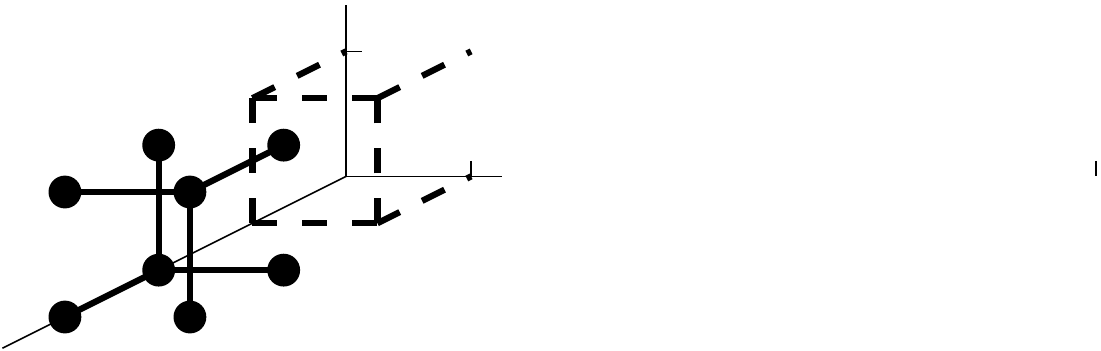}}
\hfill
\subfigure[$(2,2,1)$]{\includegraphics[scale=.4]{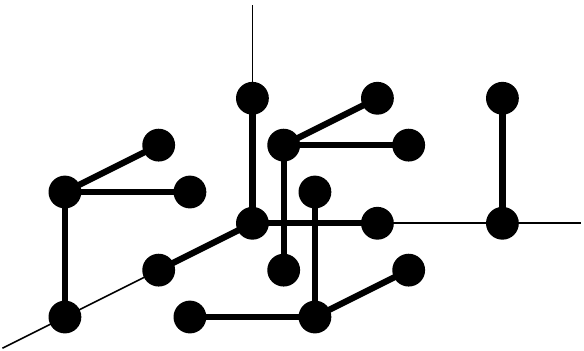}\label{fig:sl2sl2sl2_221}}
\hfill
\subfigure[$(3,2,1)$]{\includegraphics[scale=.4]{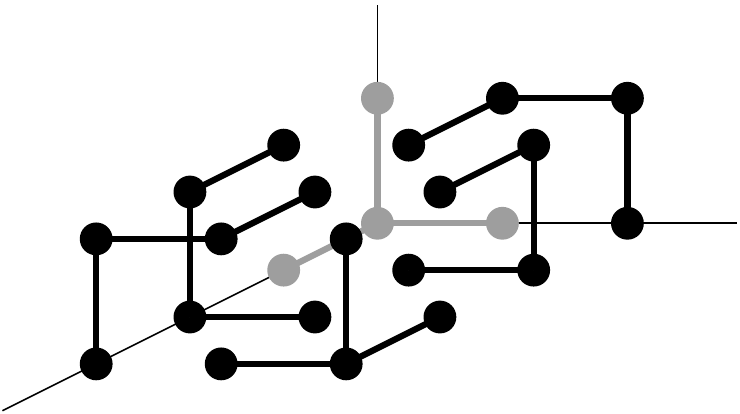}\label{fig:sl2sl2sl2_321}}
\hfill
\subfigure[$(4,2,1)$]{\includegraphics[scale=.4]{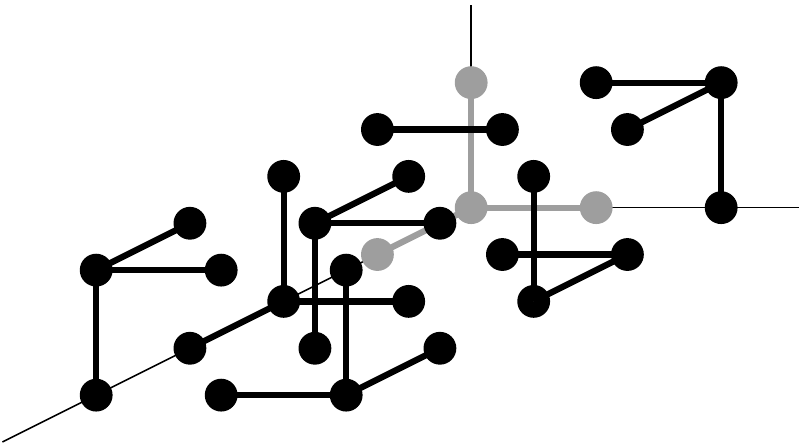}\label{fig:sl2sl2sl2_421}}
\hfill
\subfigure[$(5,2,1)$]{\includegraphics[scale=.4]{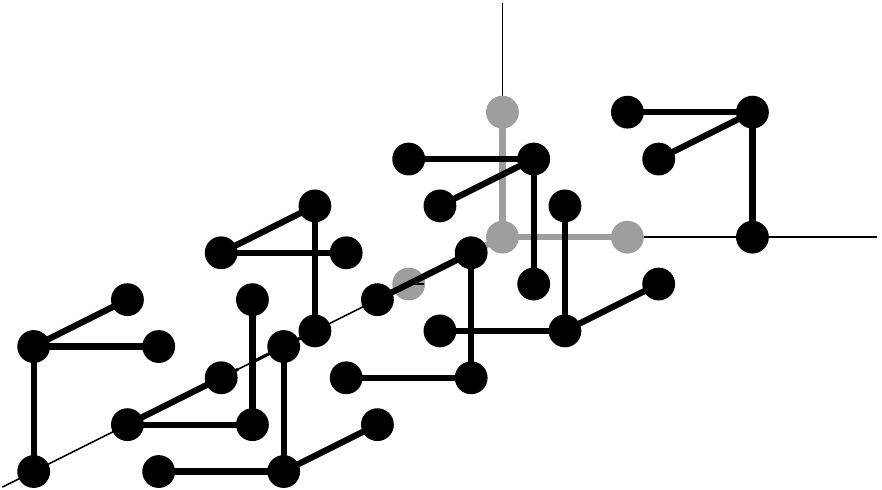}\label{fig:sl2sl2sl2_521}}
\hfill
\subfigure[induction for $(*,2,1)$]
	{\includegraphics[scale=.4]{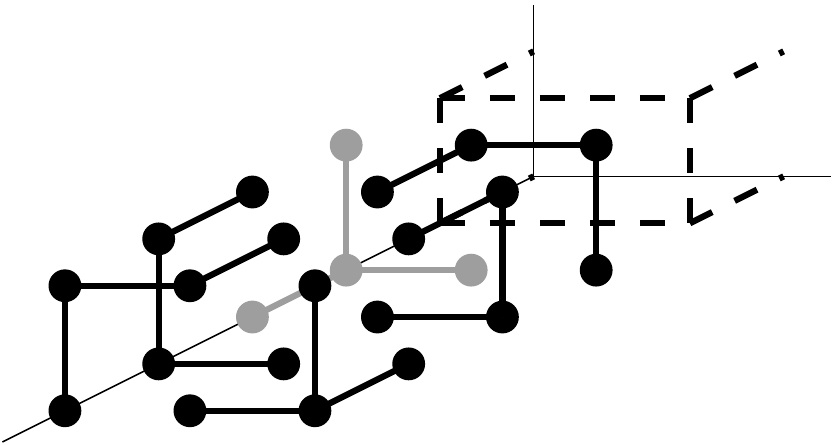}\label{fig:sl2sl2sl2_21indDEF}}
\caption{Non-defective pictures for $(e,f)=(1,1)$ or $(2,1)$}
\label{fig:sl2sl2sl2a}
\end{figure}

Figure \ref{fig:sl2sl2sl2a} gives inductive constructions for pictures
for $(e,f) \in \{(1,1),(2,1)\}$ that are non-defective except for
$(e,f)=(1,1)$ and $d$ even. The grey shades serve no other purpose than to
distinguish between front and behind. 

\begin{figure}
\subfigure[$(3,3,1)$]{\includegraphics[scale=.4]{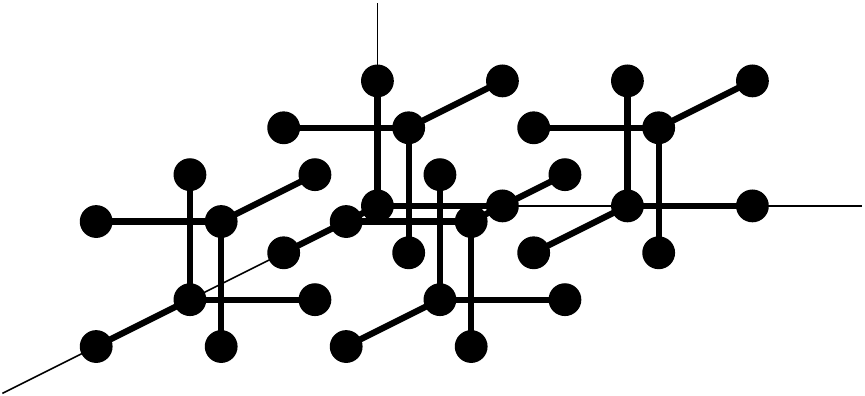}\label{fig:sl2sl2sl2_331}}
\hfill
\subfigure[induction for
$(*,3,1)$]{\includegraphics[scale=.4]{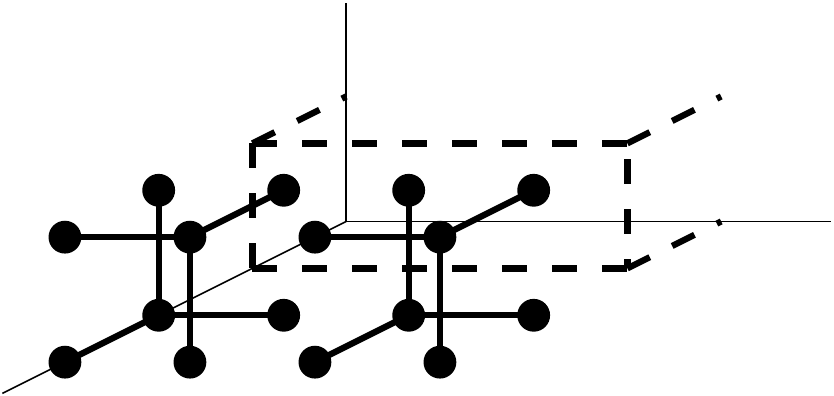}\label{fig:sl2sl2sl2_31ind}}
\hfill
\subfigure[$(4,4,1)$]{\includegraphics[scale=.4]{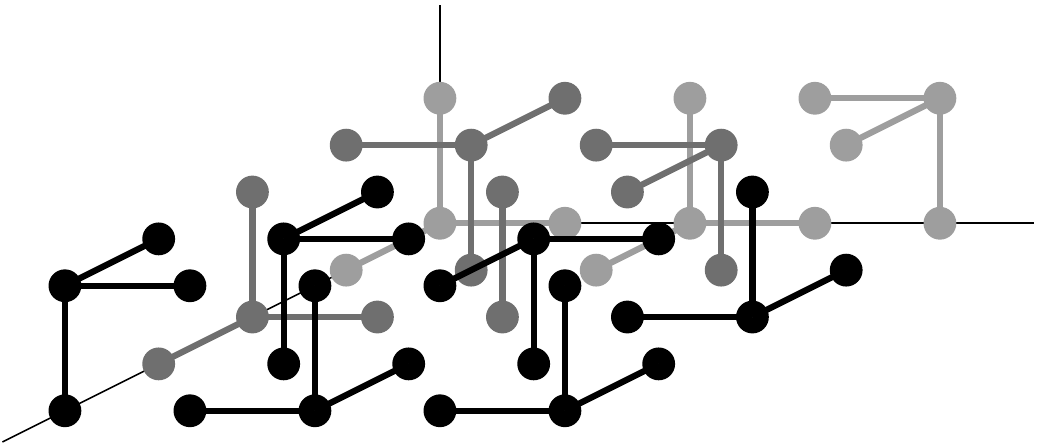}\label{fig:sl2sl2sl2_441}}
\hfill
\subfigure[induction for
$(*,4,1)$]{\includegraphics[scale=.4]{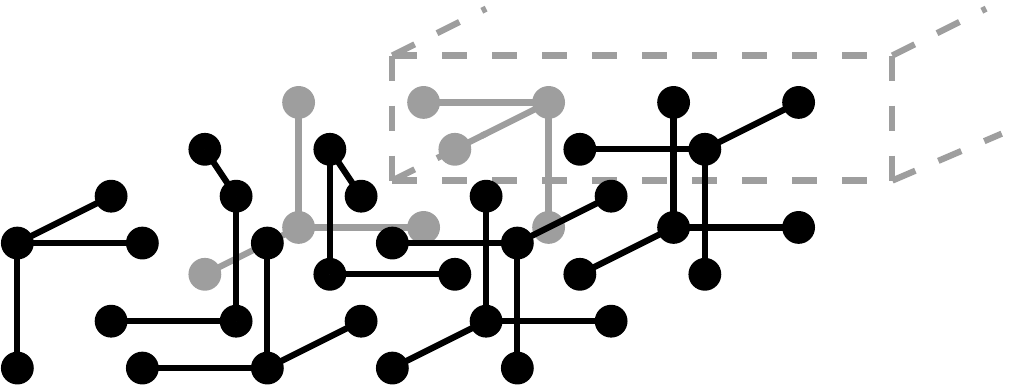}\label{fig:sl2sl2sl2_41ind}}
\caption{Non-defective pictures for $(*,3,1)$ and $(*,4,1)$.}
\label{fig:sl2sl2sl2b}
\end{figure}

Rotating appropriately, this also gives non-defective pictures for
$(1,3,1)$ and $(2,3,1)$; Figure \ref{fig:sl2sl2sl2_31ind} then gives an
inductive construction of non-defective pictures for $(d,3,1)$ for $d
\geq 3$.

So far we have found non-defective pictures for $(2,4,1)$ and
$(3,4,1)$ (just rotate those for $(4,2,1)$ and $(4,3,1)$). Figure
\ref{fig:sl2sl2sl2_441} gives a non-defective picture for $(4,4,1)$. A
non-defective picture for $(5,4,1)$ can be constructed from a
non-defective picture for $(5,1,1)$ and one for $(5,2,1)$.  Now let $d
\geq 6$ and write $d+1=4q+r$ with $q \geq 0$ and $r \in \{3,4,5,6\}$.
Then using $q$ copies of our non-defective picture for $(3,4,1)$ and
$1$ copy of our non-defective picture for $(r-1,4,1)$, we can build a
non-defective picture for $(d,4,1)$; see Figure
\ref{fig:sl2sl2sl2_41ind} for this inductive procedure.

We already have non-defective pictures for $(1,5,1)$ and $(2,5,1)$. For
$d \geq 3$, write $d+1=q*2+r$ with $r \in \{2,3\}$. Then a non-defective
picture for $(d,5,1)$ can be constructed from $q$ copies of our
non-defective picture for $(1,5,1)$ and $1$ copy of our non-defective
picture for $(r-1,5,1)$.

Let $d \geq e \geq 6$ and write $e+1=q*4+r$ with $r \in \{3,4,5,6\}$. Then
we can construct a non-defective picture for $(d,e,1)$ by putting together
$q$ non-defective pictures for $(d,3,1)$ and $1$ non-defective picture
for $(d,r-1,1)$. This settles all cases of the form $(d,e,1)$.

\begin{figure}
\subfigure[$(2,2,2)$]{\includegraphics[scale=.4]{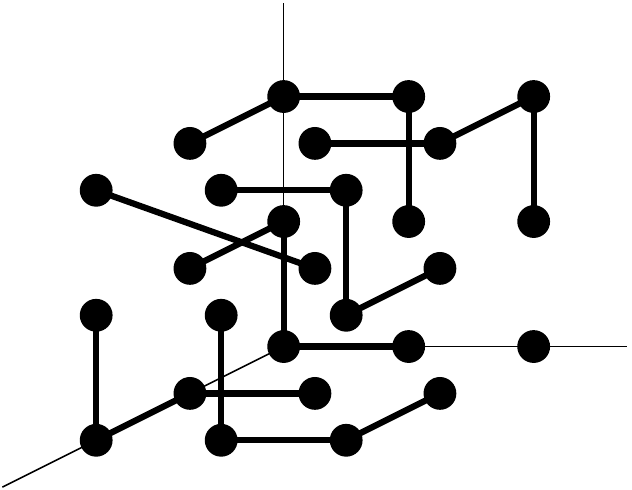}\label{fig:sl2sl2sl2_222}}
\subfigure[$(3,2,2)$]{\includegraphics[scale=.4]{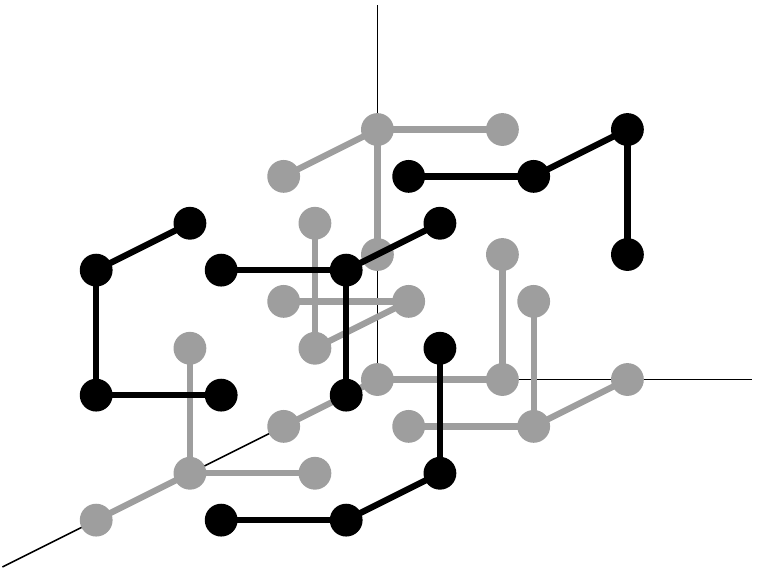}\label{fig:sl2sl2sl2_322}}
\subfigure[$(4,2,2)$]{\includegraphics[scale=.4]{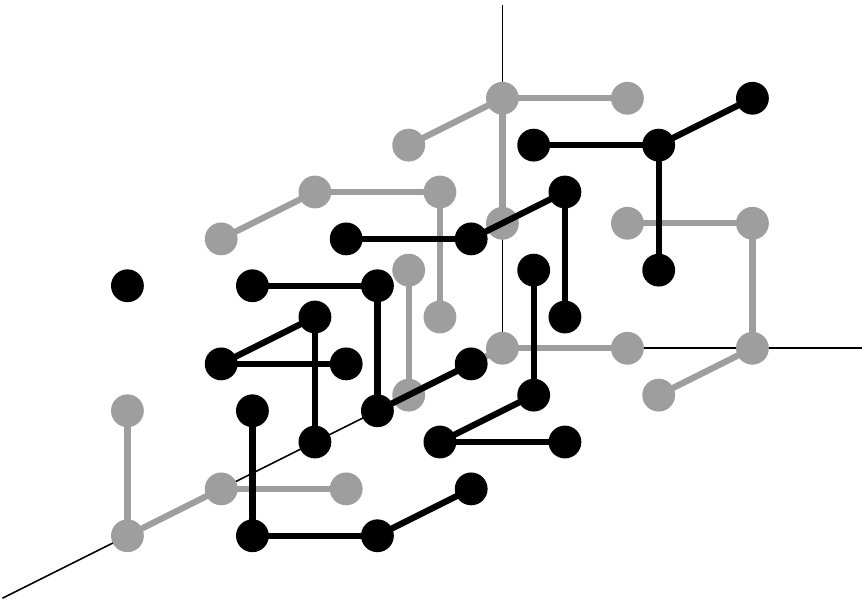}\label{fig:sl2sl2sl2_422}}
\subfigure[$(4,4,2)$]{\includegraphics[scale=.3]{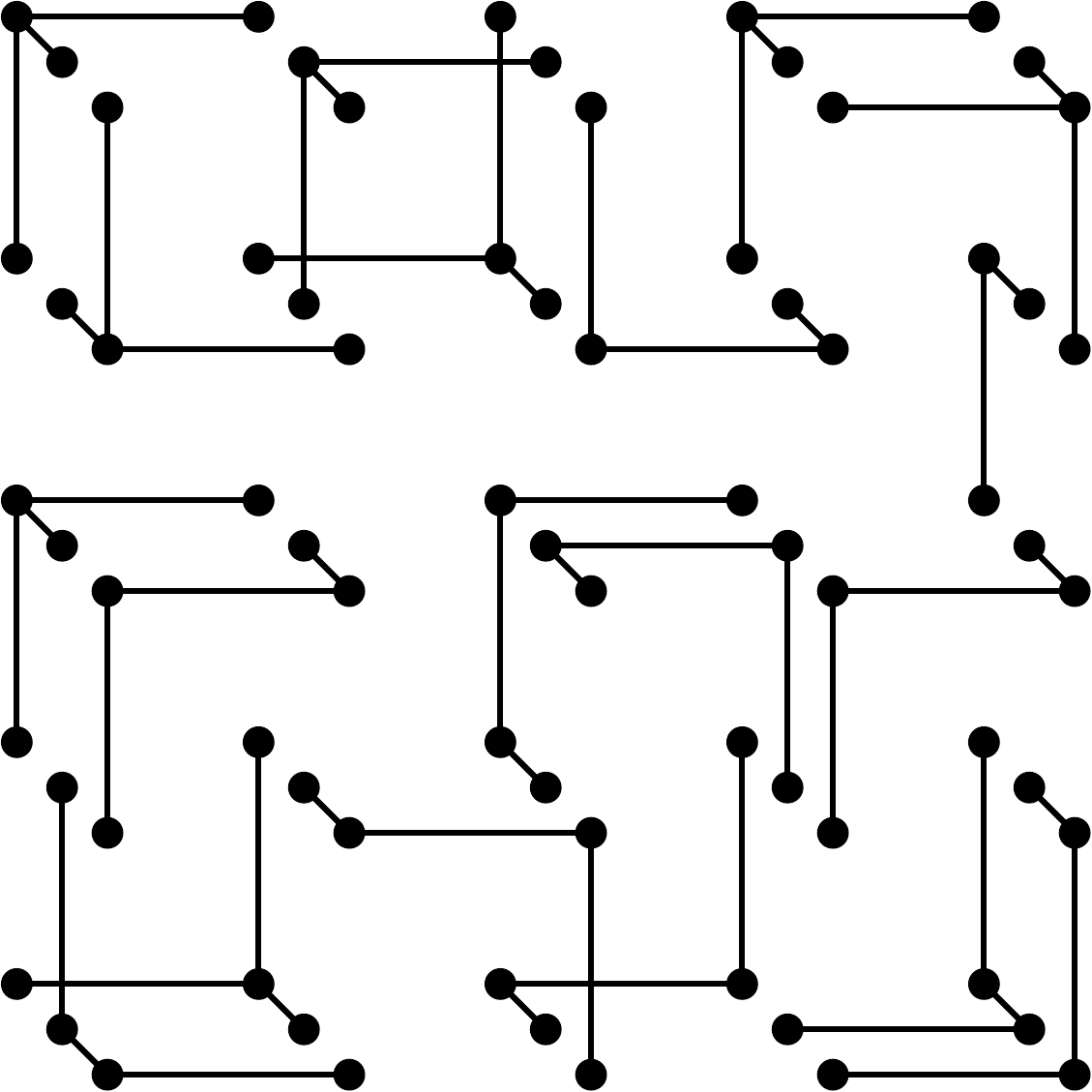}\label{fig:sl2sl2sl2_442}}
\caption{Non-defective pictures for some $(*,*,2)$.}
\label{fig:sl2sl2sl2c}
\end{figure}

Figure \ref{fig:sl2sl2sl2_222} gives a nice picture for $(d,e,f)=(2,2,2)$,
The picture is defective, but it shows that $kX$ has the expected
dimension for $k=1,\ldots,6$ and defect at most $1$ for $k=7$. From
\cite{Catalisano07} we know that $7X$ is, indeed, defective, so we are
done.  Figure \ref{fig:sl2sl2sl2_322} gives a non-defective picture for
$(3,2,2)$. Similarly, Figure \ref{fig:sl2sl2sl2_422} gives a non-defective
picture for $(4,2,2)$.

Now let $d \geq 5$ and write $d+1=(3+1)q+(r+1)$ with $r \in
\{1,3,4,6\}$. Then we can construct a non-defective picture for $(d,2,2)$
from $q$ non-defective pictures for $(3,2,2)$ and one non-defective
picture for $(r,2,2)$. This settles $(d,2,2)$. 

For $(2,3,2)$ and $(1,3,2)$ we have already found non-defective
pictures. For $d \geq 3$ write $d+1=2q + (r+1)$ with $r \in \{1,2\}$. Then
one can construct a non-defective picture for $(d,3,2)$ from $q$
non-defective pictures for $(1,3,2)$ and one non-defective picture for
$(r,3,2)$. This settles $(d,3,2)$.

If $d+1$ is even, then we can a construct non-defective picture for
$(d,e,2)$ with $d \geq e \geq 2$ as follows: write $e+1=2q + (r+1)$
with $r \in \{1,2\}$, and put together $q$ non-defective pictures for
$(d,1,2)$ and one non-defective picture for $(d,r,2)$.

Figure \ref{fig:sl2sl2sl2_442} shows how a copy of our earlier
non-defective picture for $(2,4,2)$ and a non-defective picture
for $(1,4,2)$ can be put together to a non-defective picture for
$(4,4,2)$. Now let $d \geq 6$ be even and write $d+1=4q+(r+1)$ with
$r \in \{2,4\}$. Then one can construct a non-defective picture for
$(d,4,2)$ from $q$ copies of our non-defective picture for $(3,4,2)$
and one non-defective picture for $(r,4,2)$.

This settles $(d,4,2)$.

Now suppose that $d \geq e \geq 5$ and $f=2$. Write $e+1=4*q+(r+1)$
with $r \in \{1,2,3,4\}$. Then we can construct a non-defective picture
for $(d,e,2)$ from $q$ non-defective pictures for $(d,3,2)$ and one
non-defective picture for $(d,r,2)$. This concludes the case where $d
\geq e \geq f=2$.

Consider the case where $d \geq e \geq f=3$. This case is easy now: write,
for instance, $e+1=q*2+(r+1)$ with $r \in \{1,2\}$. Then a non-defective
picture for $(d,e,3)$ can be constructed from $q$ non-defective pictures
for $(d,1,3)$ and one non-defective picture for $(d,r,3)$.

\begin{figure}
\includegraphics[scale=.4]{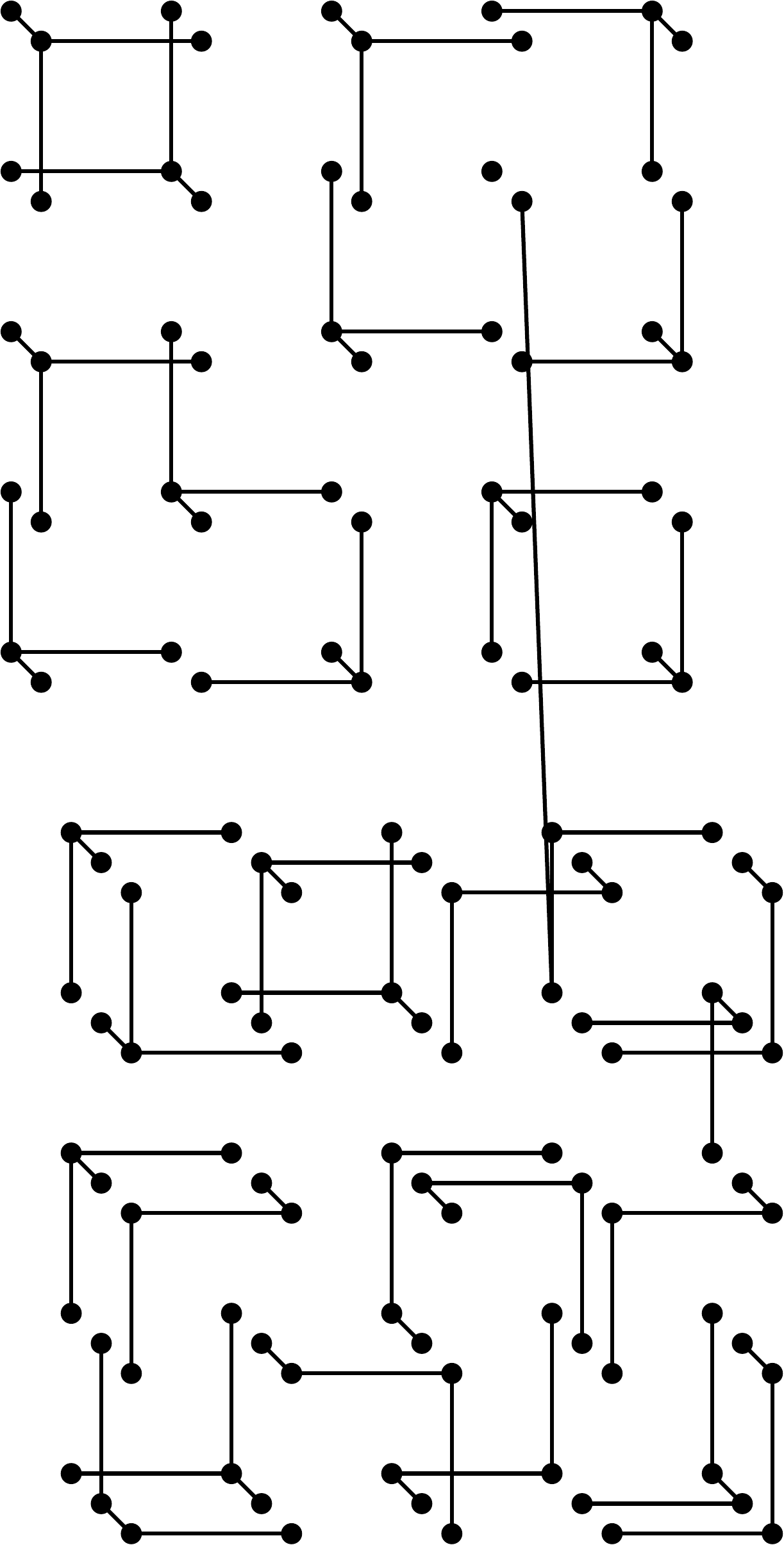}
\caption{A non-defective picture for $(4,4,4)$.}
\label{fig:sl2sl2sl2_444}
\end{figure}

The above gives (by rotating) non-defective pictures for $(d,e,4)$ for
all $d \geq 1$ and $e \in \{1,2,3\}$. Figure \ref{fig:sl2sl2sl2_444} shows
how to construct a non-defective picture for $(4,4,4)$. It may need a bit
of explanation: the upper half is a non-defective picture for $(4,4,1)$,
very close to that in Figure \ref{fig:sl2sl2sl2_441}---but the superflous
pair of vertices is separated. The lower half is a non-defective picture
for $(4,4,2)$, very close to that in Figure \ref{fig:sl2sl2sl2_442}. By
joining the lower one of the superflous vertices in the upper half with
the triangle in the lower half, we create a non-defective picture for
$(4,4,4)$. Now suppose that $d \geq e \geq 5$ and write $e+1=4q+(r+1)$
with $r \in \{1,2,3,4\}$. Then we find a non-defective picture for
$(d,e,4)$ from $q$ non-defective pictures for $(d,3,4)$ and one
non-defective picture for $(d,r,4)$.

Finally, suppose that $d \geq e \geq f \geq 5$, and write $f+1=4q+(r+1)$
with $r \in \{1,2,3,4\}$. Then a non-defective picture for $(d,e,f)$
can be assembled from $q$ non-defective pictures for $(d,e,3)$ and one
non-defective picture for $(d,e,r)$. This concludes the
proof of Theorem \ref{thm:P2P1}.

\section{Secant dimensions of $\PP^2 \times \PP^1$}
\label{sec:P2P1}

For Theorem~\ref{thm:P2P1} we first deal with the defective cases: the
Segre-Veronese embeddings of degree $(2,\text{even})$ are all defective
by \cite[Example 3.2]{Catalisano06}. That the embedding of degree $(3,1)$
is defective can be proved using a polynomial interpolation argument,
used in \cite{Catalisano05c} for proving defectiveness of other secant
varieties: Split $(3,1)=(2,0)+(1,1)$. Now it is easy to see that given
$5$ general points there exist non-zero forms $f_1,f_2$ of multi-degrees
$(2,0)$ and $(1,1)$, respectively, that vanish on those points. But
then the product $f_1f_2$ vanishes on those points together with all its
first-order derivatives; hence the $5$-th secant variety does not fill the
space. The proof below shows that its codimension is not more than $1$.

For the non-defective proofs we have to solve the optimisation problems
$\AP(B,k)$, where
\[ B=\{(x,y,z) \in \ZZ^3 \mid x,y,z \geq 0, x+y \leq d,
	\text{ and } z \leq e\}. \]
We will do a double induction over the degrees $e$ and $d$: First, in
Subsections \ref{ssec:e1}---\ref{ssec:e4} we treat the cases where $e$
is fixed to $1,2,3,4$, respectively, by induction over $d$. Then, in
Subsection \ref{ssec:d} we perform the induction over $e$. We will always
think of the $x$-axis as pointing towards the reader, the $y$-axis as
pointing to the right and the $z$-axis as the vertical axis. By $T_{d,e}$
we will mean a picture (non-defective, if possible) for $(d,e)$. We will
also use (non-defective) pictures from Section \ref{sec:P1P1P1} as building
blocks; we denote the picture for the Segre-Veronese embedding of $\PP^1
\times \PP^1 \times \PP^1$ by of degree $(a,b,c)$ by $B_{a,b,c}$.

\subsection{The cases where $e=1$} \label{ssec:e1}

\begin{figure}
\subfigure[$T_{1,1}$] {\includegraphics{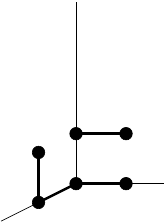}\label{fig:101}}
\subfigure[$T_{2,1}$] 
	{\includegraphics[scale=.5]{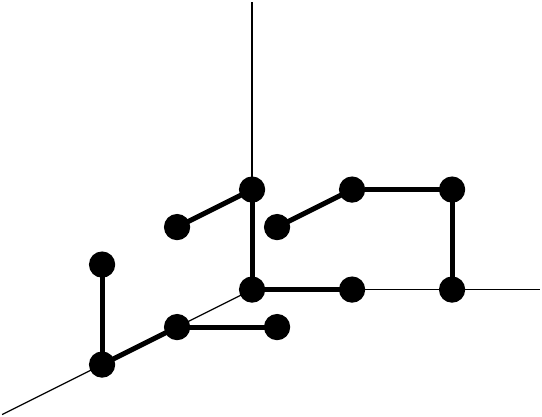}\label{fig:201}}
\subfigure[$T_{3,1}$; defective]
      {\includegraphics[scale=.8]{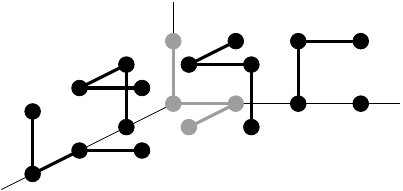}\label{fig:301}}
\subfigure[$T_{4,1}$]
      {\includegraphics[scale=.5]{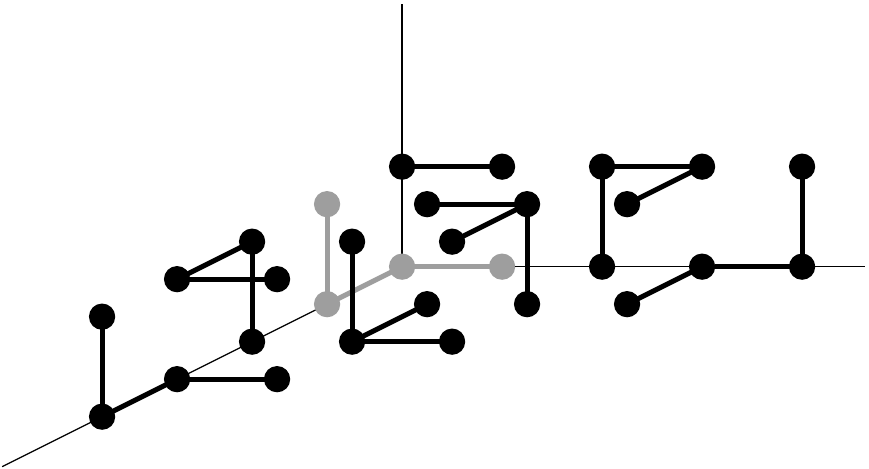}\label{fig:401}}
\caption{Induction basis for $(*,1)$} 
\end{figure}

\begin{figure}[ht]
\begin{center}
\includegraphics{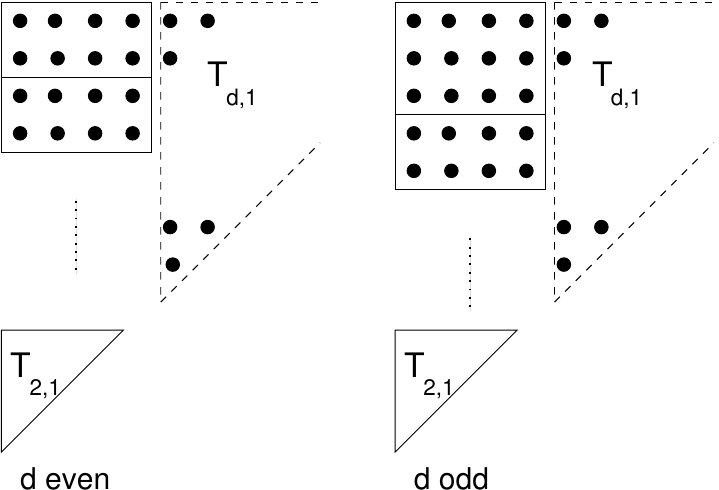}
\caption{Induction step for $(*,1)$}\label{fig:induction-1}
\end{center}
\end{figure}

Figures \ref{fig:101}--\ref{fig:401} give pictures for $(d,1)$ with
$d=1,\ldots,4$. Now we explain how to construct a non-defective picture
for $(d+4,1)$ from a non-defective picture for $(d,1)$: First translate
$T_{d,1}$ four steps to the right, and then proceed as follows:
\begin{enumerate}
\item If $d+1$ is even, $d+1=2l$ for some $l$, then put $l+1$ copies
of $B_{1,3,1}$ to the left of $T_{d,1}$, starting at the origin. Finally,
add a copy of $T_{2,1}$. 
\item If $d+1$ is odd, $d+1=2l+1$ for some $l$, then put one copy of
$B_{2,3,1}$ and $l-1$ copies of $B_{1,3,1}$ to the left of $T_{d,1}$.
Finally, add another copy of $T_{2,1}$.
\end{enumerate}
This is illustrated in Figure~\ref{fig:induction-1}.

To complete the induction, since $T_{3,1}$ is defective, we need a
non-defective picture for $(7,1)$. We can construct this by using
two copies of $T_{2,1}$, a box $B_{1,3,1}$ and the block $B_{2,2,1}$
(at the origin) from Figure~\ref{fig:sl2sl2sl2_221}.  The remaining
vertices are grouped together as in Figure~\ref{fig:701} below.

\begin{figure}[ht]
\begin{center}
\includegraphics{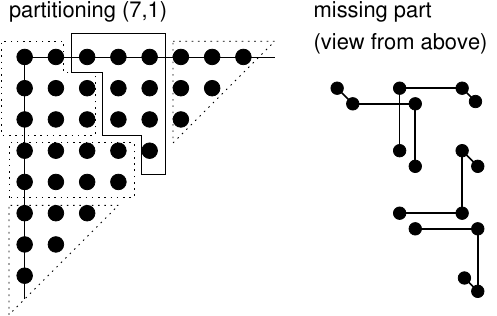}
\caption{Obtaining $T_{7,1}$}\label{fig:701}
\end{center}
\end{figure}

\subsection{The cases where $e=2$} \label{ssec:e2}

\begin{figure}
\begin{center}
\subfigure[$T_{1,2}$]
	{\includegraphics[scale=1]{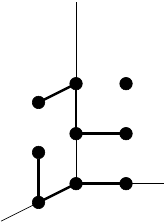}\label{fig:102}}
\subfigure[$T_{2,2}$; defective]
	{\includegraphics[scale=.5]{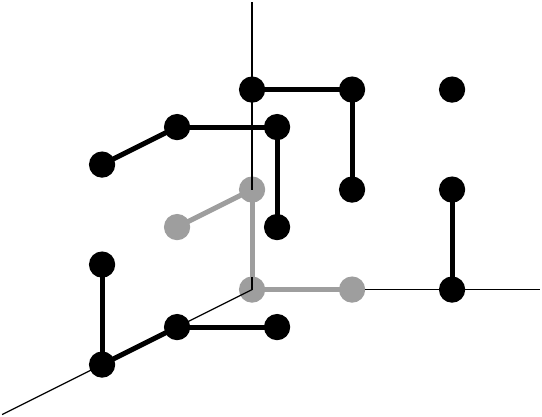}\label{fig:202}}
\subfigure[$T_{3,2}$]
	{\includegraphics[scale=.8]{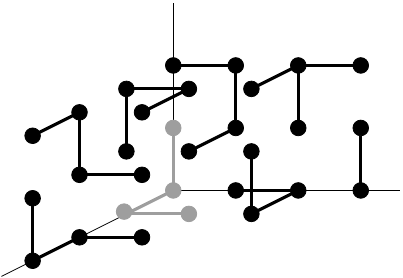}\label{fig:302}}
\subfigure[$T_{4,2}$]
	{\includegraphics{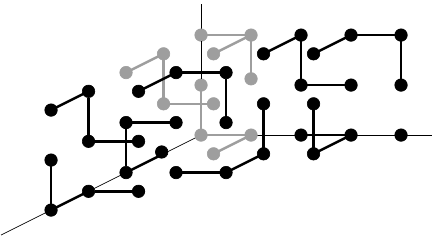}\label{fig:402}}
\subfigure[$T_{5,2}$]
     	{\includegraphics{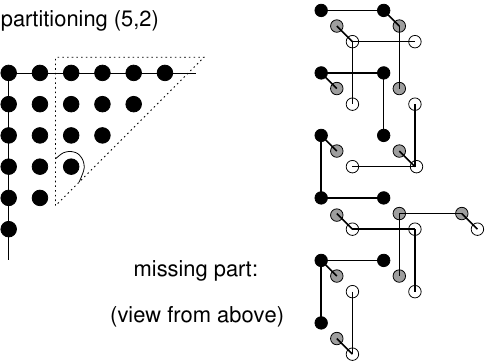}\label{fig:502}}
	\hfill
\subfigure[$T_{6,2}$]
     	{\includegraphics{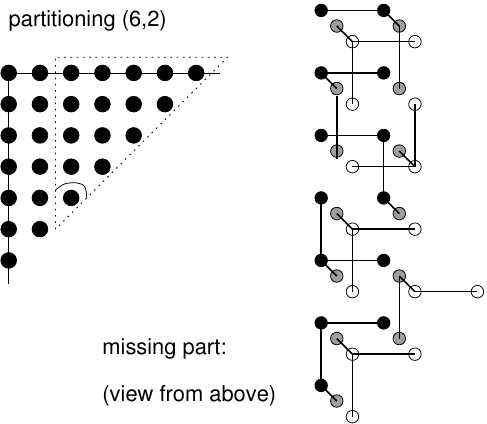}\label{fig:602}}
\subfigure[$T_{7,2}$]
     	{\includegraphics{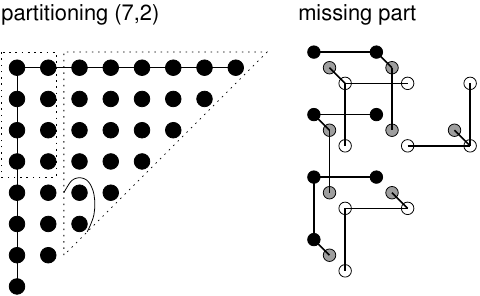}\label{fig:702}}
	\hfill
\subfigure[$T_{8,2}$]
     	{\includegraphics{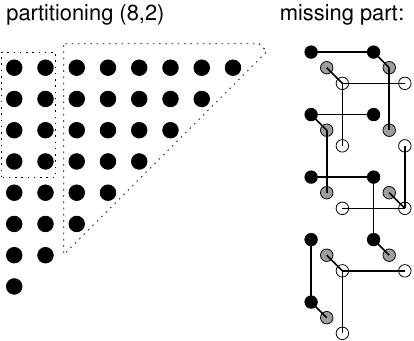}\label{fig:802}}
\subfigure[$T_{10,2}$]
     	{\includegraphics{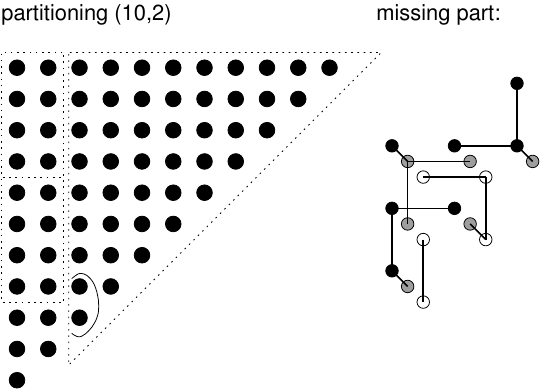}\label{fig:1002}}
\end{center}
\caption{Induction basis for $(*,2)$.}
\end{figure}

Figures \ref{fig:102}-\ref{fig:1002} lay the basis for the induction over
$d$. Note that $T_{6,2}$ is the first among the pictures whose number
of vertices is divisible by $4$. To finish the induction, we need to
construct a non-defective picture for $(d+8,2)$ from $T_{d,2}$. First of
all, move $T_{d,2}$ eight positions to the right. Then proceed as follows:
\begin{enumerate} 
\item If $d$ is odd, $d=2l+1$ for some $l\ge 0$, put $l$ pairs
of $B_{1,3,2}$ to the right of $T_{d,2}$ (starting at the origin),
then two copies of $B_{2,3,2}$, and finally a copy of $T_{6,2}$.

\item If $d$ is even, $d=2l$ for some $l>0$, put $l+1$ pairs
of $B_{1,3,2}$ starting at the origin.  Finish off with one copy of
$T_{6,2}$. 
\end{enumerate}

\begin{figure}[ht]
\begin{center}
\includegraphics{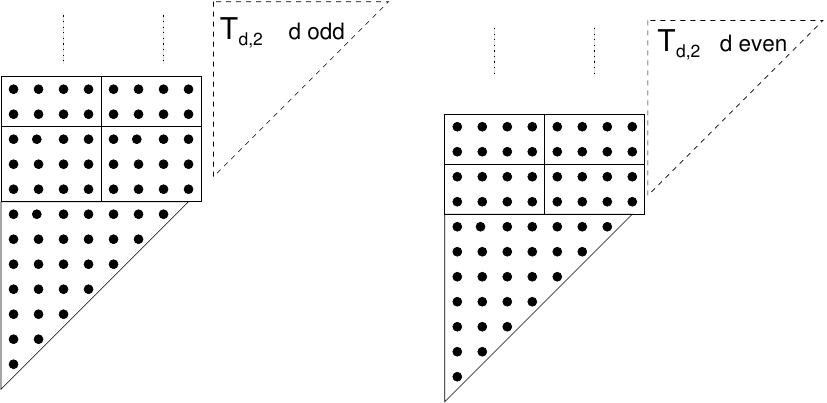}
\caption{Induction step for $(*,2)$}\label{fig:induction-2}
\end{center}
\end{figure}

This is illustrated in Figure~\ref{fig:induction-2}.

\subsection{The cases where $e=3$} \label{ssec:e3}

\begin{figure}
\begin{center}
\subfigure[$T_{1,3}$]
      {\includegraphics[scale=.8]{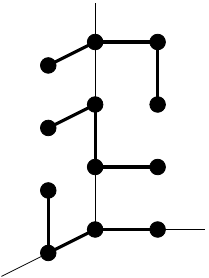}\label{fig:103}}
\subfigure[$T_{2,3}$]
      {\includegraphics[scale=.5]{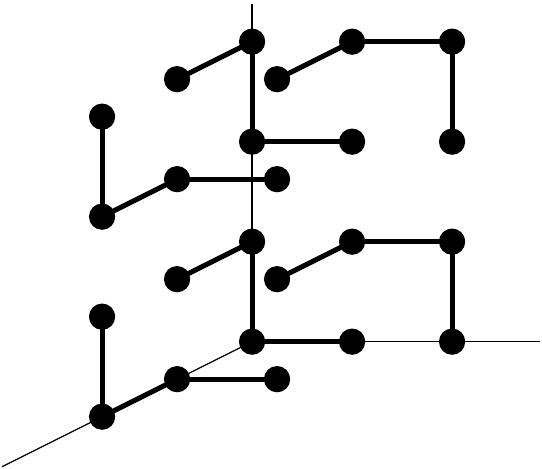}\label{fig:203}}
\caption{Induction basis for $(*,3)$}
\end{center}
\end{figure}

Here the induction over $d$ is easier since every $T_{d,3}$ has its number
of vertices divisible by $4$. Figures \ref{fig:103} and \ref{fig:203}
lay the basis of the induction (the latter just consists of two copies
of $T_{2,1}$). Now we show that from a non-defective
$T_{d,3}$ with $d$ odd one can construct non-defective
$T_{d+2,3}$ and $T_{d+3,3}$. Write $d=2l+1$, and proceed as
follows. 
\begin{enumerate}
\item Move $T_{2l+1,3}$ two positions to the right. Put a block
$B_{2l+1,1,3}$ at the origin, and conclude with a copy of $T_{1,3}$. This
gives $T_{2l+3,3}$.

\item Move $T_{2l+1,3}$ three steps to the right. Put a block
$B_{2l+1,2,3}$ at the origin, and conclude with a copy of $T_{2,3}$.
\end{enumerate}

\begin{figure}[ht]
\begin{center}
\subfigure[odd $d+2$]
      {\includegraphics{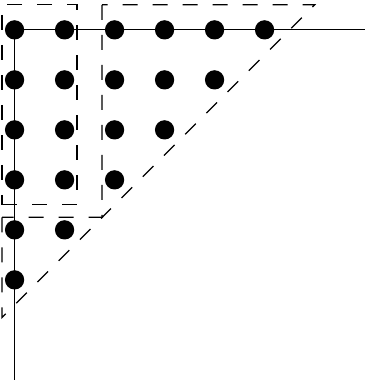}\label{fig:add-odd}}
   \hspace{1cm}
\subfigure[even $d+3$]
      {\includegraphics{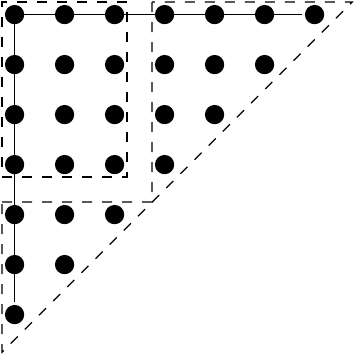}\label{fig:add-even}}
\caption{Induction steps from $T_{3,3}$}\label{fig:to-3_3}
\end{center}
\end{figure}

For $d=3$ this is illustrated in Figure~\ref{fig:to-3_3}.

\subsection{The cases where $e=4$} \label{ssec:e4}

\begin{figure}
\begin{center}
\subfigure[$T_{1,4}$]
	{\includegraphics{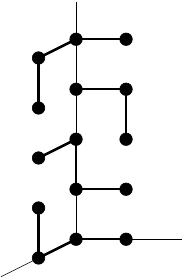}\label{fig:104}}
\subfigure[$T_{2,4}$; defective]
	{\includegraphics[scale=.5]{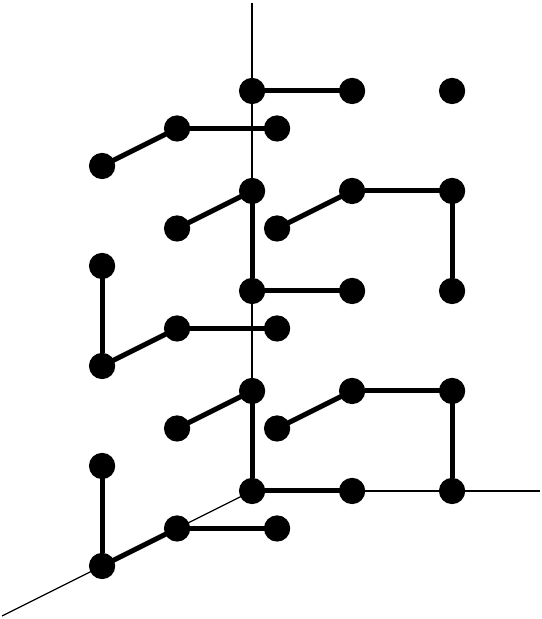}\label{fig:204}}
\subfigure[$T_{3,4}$]
	{\includegraphics[scale=.5]{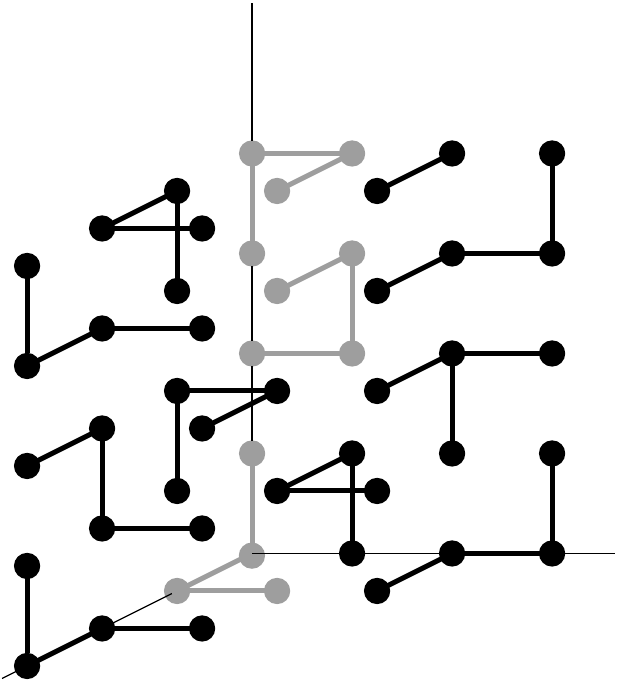}\label{fig:304}}
\subfigure[$T_{4,4}$]
	{\includegraphics{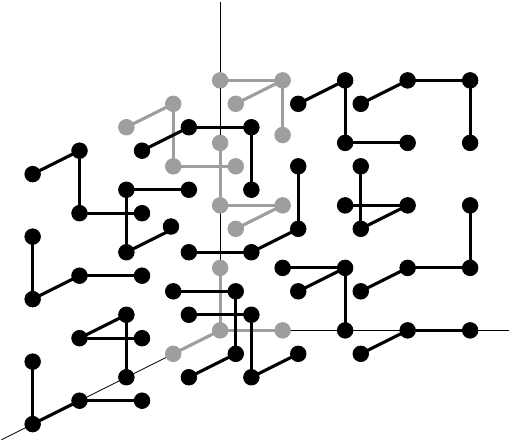}\label{fig:404}}
	\hfill
\subfigure[$T_{5,4}$]
	{\includegraphics{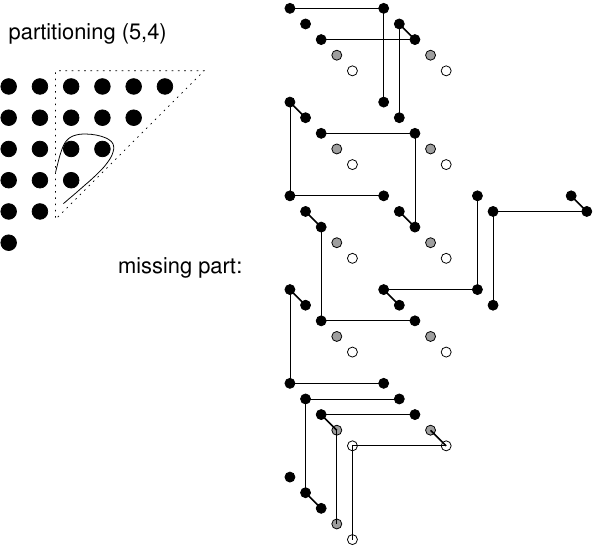}\label{fig:504}}
\subfigure[$T_{6,4}$]
	{\includegraphics{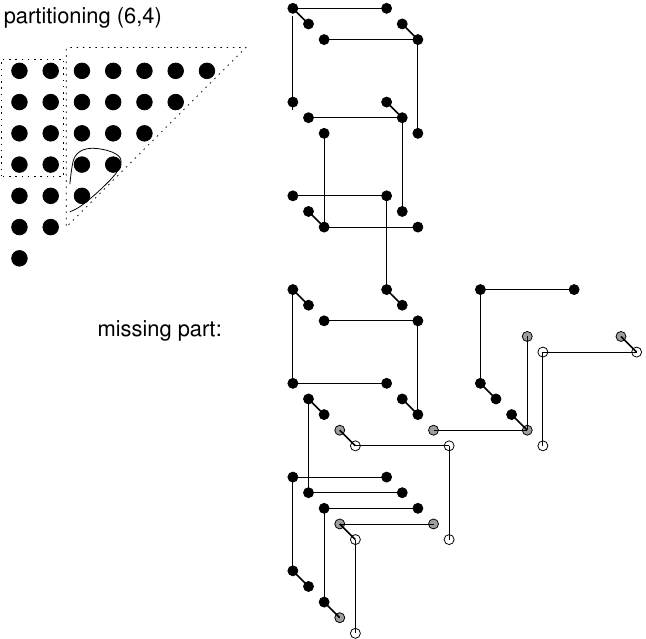}\label{fig:604}}
	\hfill
\subfigure[$T_{7,4}$]
	{\includegraphics{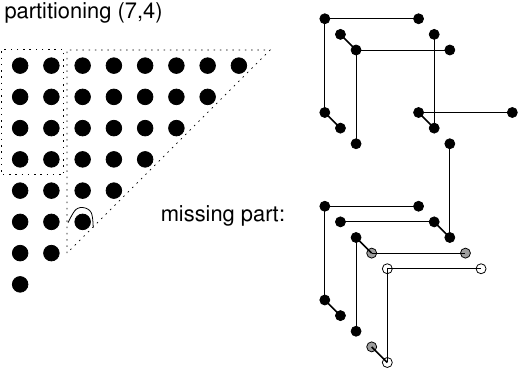}\label{fig:704}}
\subfigure[$T_{8,4}$]
	{\includegraphics{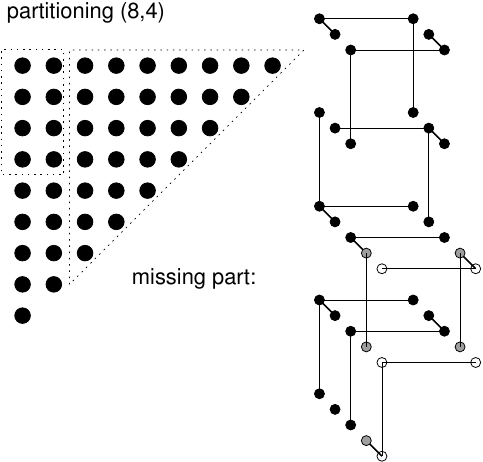}\label{fig:804}}
	\hfill
\subfigure[$T_{10,4}$]
	{\includegraphics{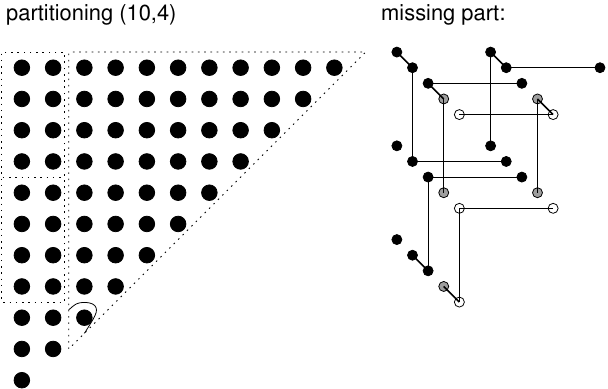}\label{fig:1004}}
\end{center}
\caption{Induction basis for $(*,4)$}
\end{figure}

Figures \ref{fig:104}-\ref{fig:1004} lay the basis of the induction. The
induction step is identical to that where $e=2$, except that the blocks
$B_{1,3,2}$ and $B_{2,3,2}$ have to be replaced by the blocks $B_{1,3,4}$
and $B_{2,3,4}$, and $T_{6,2}$ has to be replaced by $T_{6,4}$.

\subsection{Induction over $e$} \label{ssec:d}

\begin{figure}
\begin{center}
\includegraphics{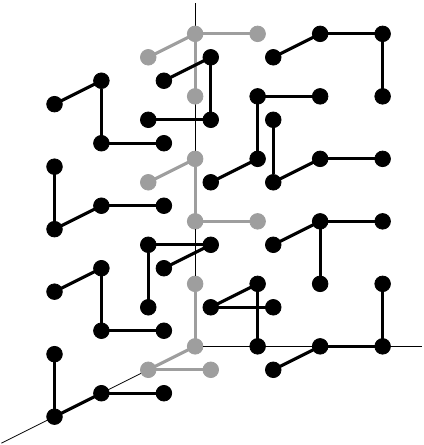}
\caption{$T_{3,5}$}\label{fig:305}
\end{center}
\end{figure}

From a non-defective picture for $(d,e)$ we can construct a non-defective
picture for $(d,e+4)$ by stacking a non-defective picture for $(d,3)$,
whose number of vertices is divisible by $4$, on top of it.  This settles
all $(d,e)$ except for those that are modulo $(0,4)$ equal to the
defective $(3,1)$ or $(2,2)$. The latter are easily handled, though:
stacking copies of $T_{2,1}$ on top of $T_{2,2}$ gives pictures for all
$(2,e)$ with $e$ even that are defective but give the correct, known,
secant dimensions. So to finish our proof of Theorem \ref{thm:P2P1} we
only need the non-defective picture for $(3,5)$ of Figure \ref{fig:305}.

\section{Secant dimensions of the point-line flag variety $\mF$} \label{sec:F}

In this section, $X=\mF$, $G=\lieg{SL}_3$, and the highest weight
$\lambda$ equals $m \omega_1 + n \omega_2$ with $m,n>0$. We first
argue that $(m,n)=(1,1)$ and $(m,n)=(2,2)$ yield defective embeddings
of $\mF$. The first weight is the adjoint weight, so the cone $C_{1,1}$
over the image of $\mF$ is just the set of rank-one, trace-zero matrices
in $\liea{sl}_3$, whose secant dimensions are well known. For the second
weight let $C_{2,2}$ be the image of $C_{1,1}$ under the map $\liea{sl}_3
\rightarrow S^2(\liea{sl}_3), v \mapsto v^2$. Then $C_{2,2}$ spans
the $\lieg{SL}_3$-submodule (of codimension $9$) in $S^2(\liea{sl}_3)$
of highest weight $2 \omega_1 + 2 \omega_2$, while it is contained in
the quadratic Veronese embedding of $\liea{sl}_3$. Viewing the elements
of $S^2(\liea{sl}_3)$ as symmetric $8 \times 8$-matrices, we find that
$C_{2,2}$ consists of rank $1$ matrices, while it is not hard to prove
that the module it spans contains matrices of full rank $8$.  Hence $7
C_{2,2}$ cannot fill the space.

For the non-defective proofs let $\alpha_1,\alpha_2$ be the
simple positive roots, so that $\Xu=\{\beta_1,\beta_2,\beta_3\}$
with $\beta_1=-\alpha_1$, $\beta_2=-\alpha_1 -\alpha_2$ and
$\beta_3=-\alpha_2$.  The subscripts indicate the order in which
the PBW-monomials are computed: for $r=(n_1,n_2,n_3)$ we write
$m_r:=X_{\beta_1}^{n_1} X_{\beta_2}^{n_2} X_{\beta_3}^{n_3}
v_\lambda$. Set
\[ B:=\{(n_1,n_2,n_3)  \in \ZZ^3 \mid
        0 \leq n_2 \leq m,
        0 \leq n_3 \leq n, \text{ and }
        0 \leq n_1 \leq m + n_3 - n_2 \}, \]
and let $M$ be the set of all $r \in \NN^{\Xu}$ with $m_r \neq 0$. We
will not need $M$ explicitly; it suffices to observe that $r_3 \leq n$
for all $r \in M$: indeed, if $r_3 > n$ then $X_{\beta_3}^{r_3} v_\lambda$
is already $0$, hence so is $m_r$. We use the following consequence of
the theory of canonical basis; see \cite[Example 10, Lemma 11]{Graaf03}.

\begin{lm}
The $m_b, b\in B$ form a basis of $V$.
\end{lm}

\begin{re} \label{re:Transposition}
The map $(n_1,n_2,n_3) \mapsto (n_1,n-n_3,m-n_2)$ sends the set $B$, which 
corresponds on the highest weight $(m,n)$, to the set corresponding 
to the highest weight $(n,m)$. Hence if we have a non-defective picture for 
one, then we also have a non-defective picture for the other. We will use 
this fact occasionally.
\end{re}

We want to apply Lemma \ref{lm:Relation}. First note that $r,r' \in
\RR^{\Xu}=\RR^3$ have the same weight if and only if $r-r'$ is a scalar
multiple of $z:=(1,-1,1)$. We set $r < r'$ if and only if $r-r'$ is a
{\em positive} scalar multiple of $z$.

\begin{lm} \label{lm:LowerIdeal}
For all $r \in M \setminus B$ and all $b\in B$ with $\xi(b)=\xi(r)$
we have $b<r$, i.e., the difference $b-r$ is a positive scalar multiple
of $z$.
\end{lm}

\begin{proof}
Suppose that $b=(n_1,n_2,n_3) \in B$ and that $n_3<n$. Then the defining
inequalities of $B$ show that $b+z=(n_1+1,n_2-1,n_3+1)$ also lies in $B$.
This shows that $B$ is a {\em lower ideal} in $(M,\leq)$, i.e., if $b
\in B$ and $r \in M$ with $r< b$, then also $r \in B$. This readily
implies the lemma.
\end{proof}

\begin{prop}
$\VP(B,k)$ is a lower bound on $\dim kC$ for all $k$.
\end{prop}

\begin{proof}
This follows immediately from Lemma \ref{lm:Relation} and Lemma
\ref{lm:LowerIdeal} when we take for $Z$ the one-dimensional cone
$\RR_{\geq 0} \cdot z$.
\end{proof}

\begin{figure}
\begin{center}
\subfigure[$(2,1)$]{\includegraphics[scale=.4]{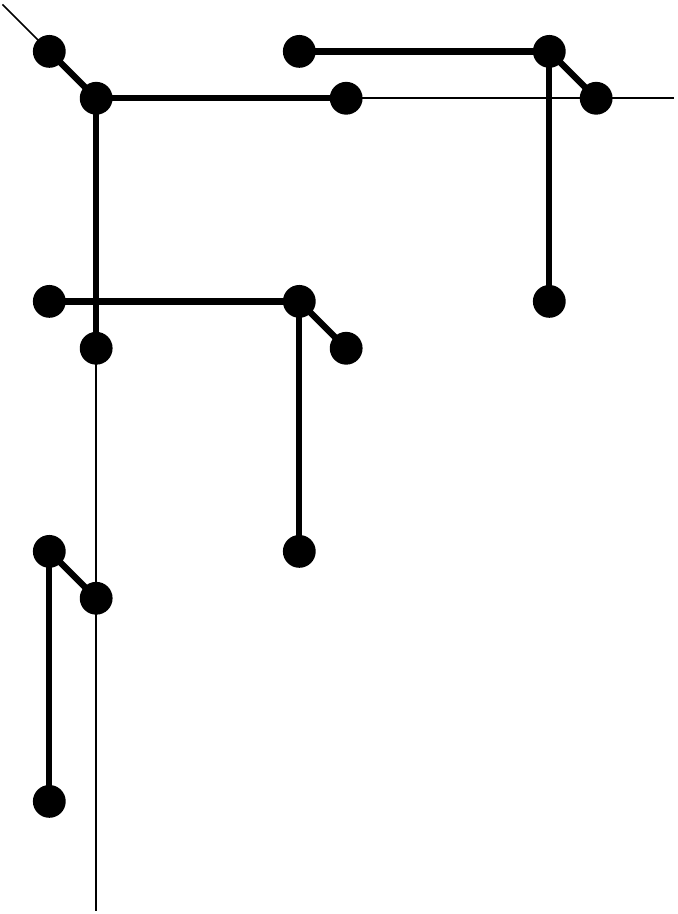}\label{fig:sl3_21}}
\subfigure[$(2+2k,1)$]{\includegraphics[scale=.4]{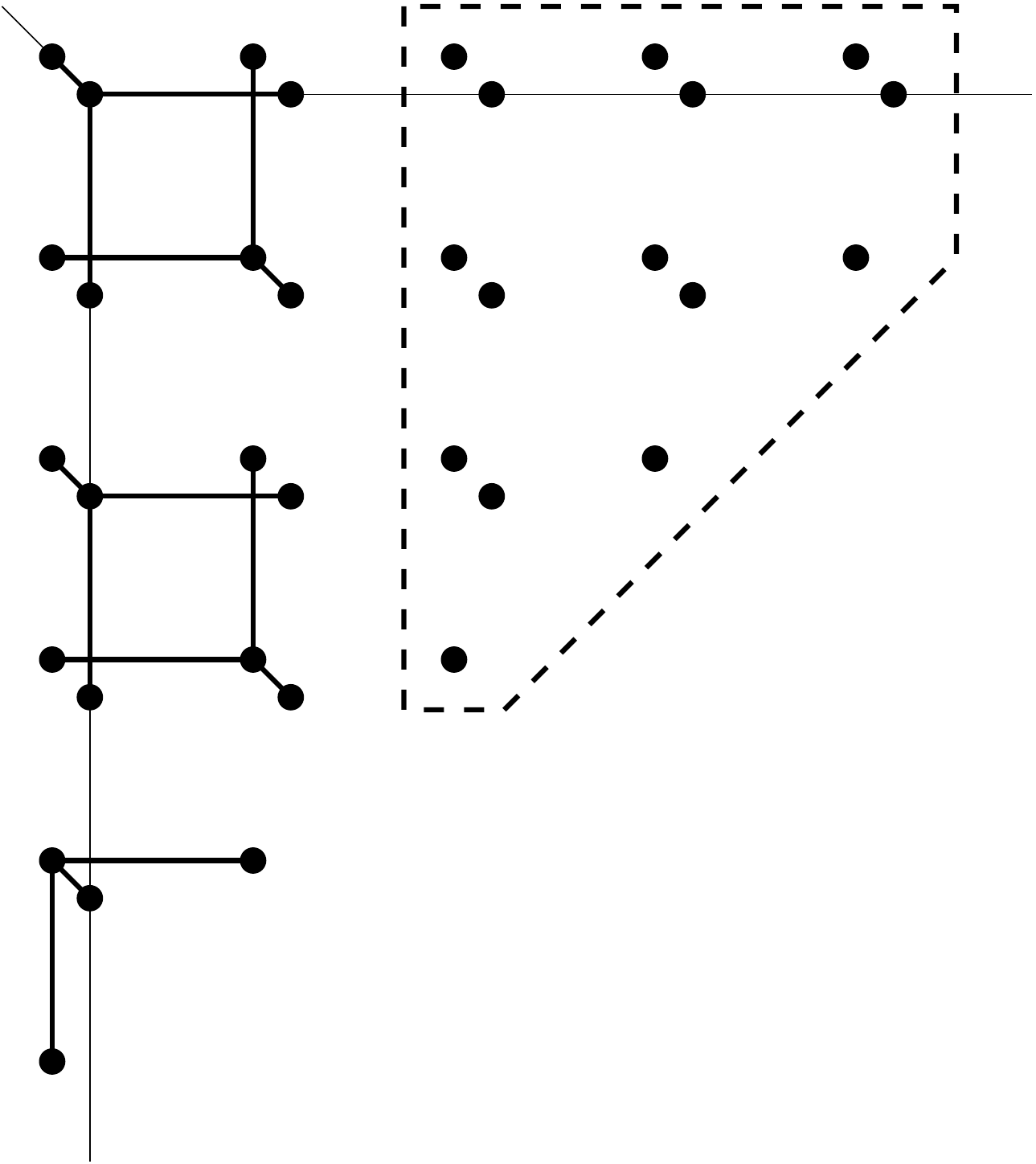}\label{fig:sl3_1ind1}}
\caption{Pictures for $(\text{even},1)$.}\label{fig:sl3_1xeven}
\end{center}
\end{figure}

\begin{figure}
\begin{center}
\subfigure[$(1,1)$;
defective]{\includegraphics[scale=.3]{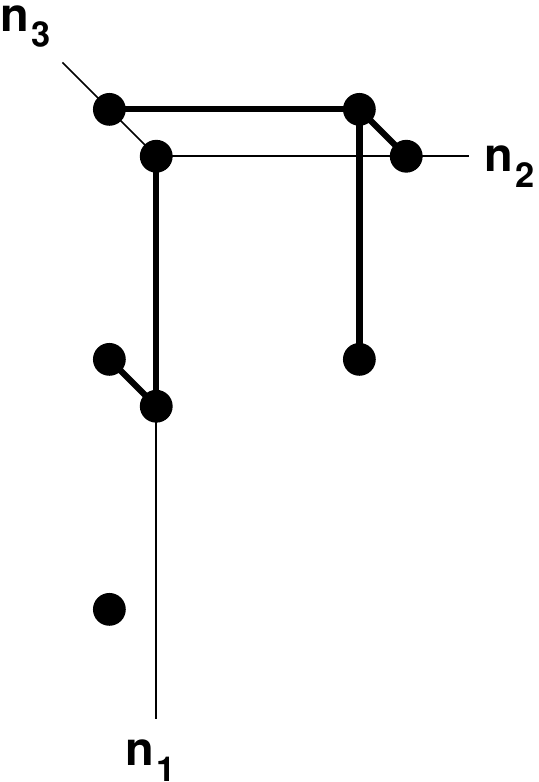}\label{fig:sl3_11}}
\subfigure[$(3,1)$]{\includegraphics[scale=.3]{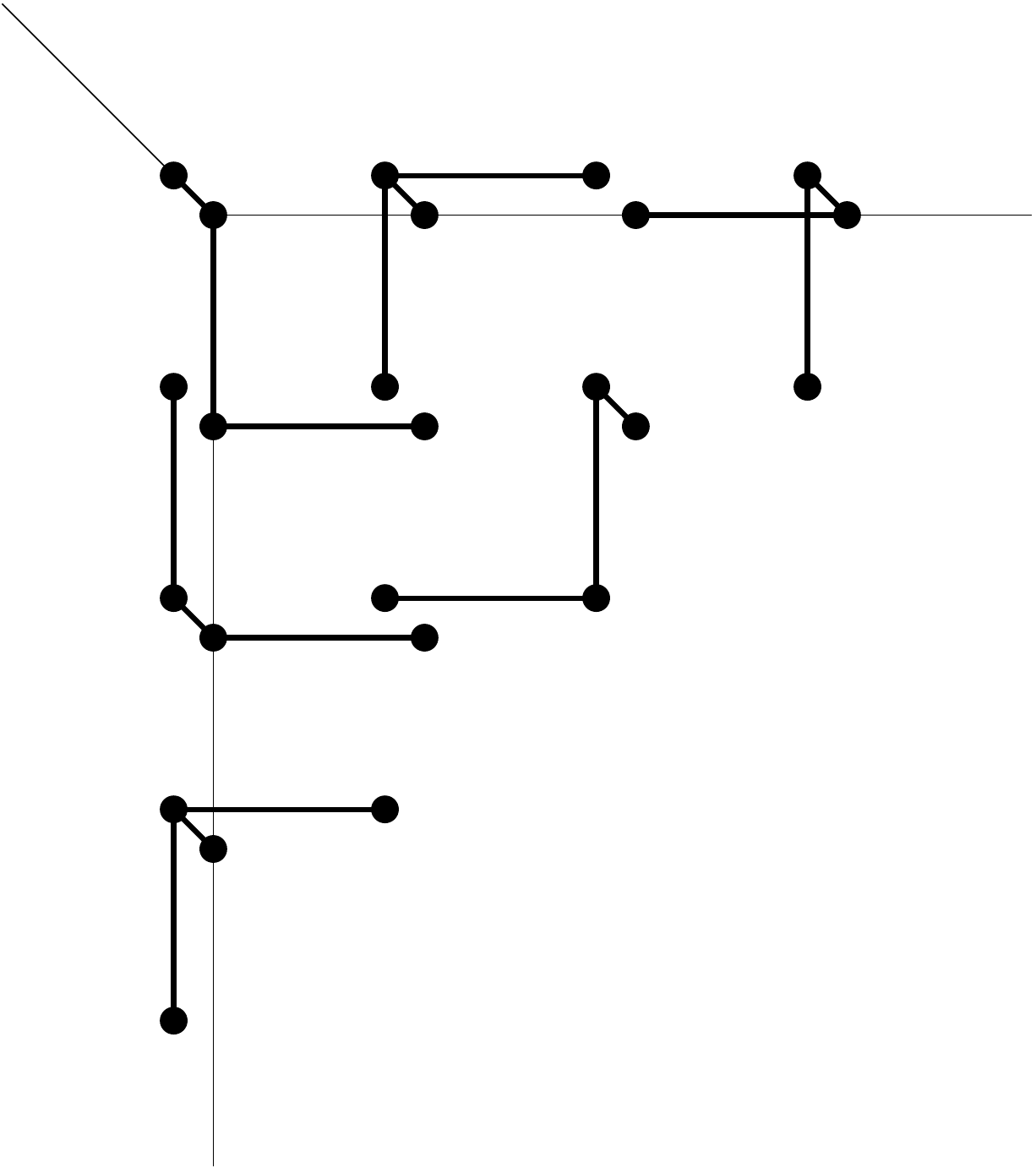}\label{fig:sl3_31}}
\subfigure[$(5,1)$]{\includegraphics[scale=.3]{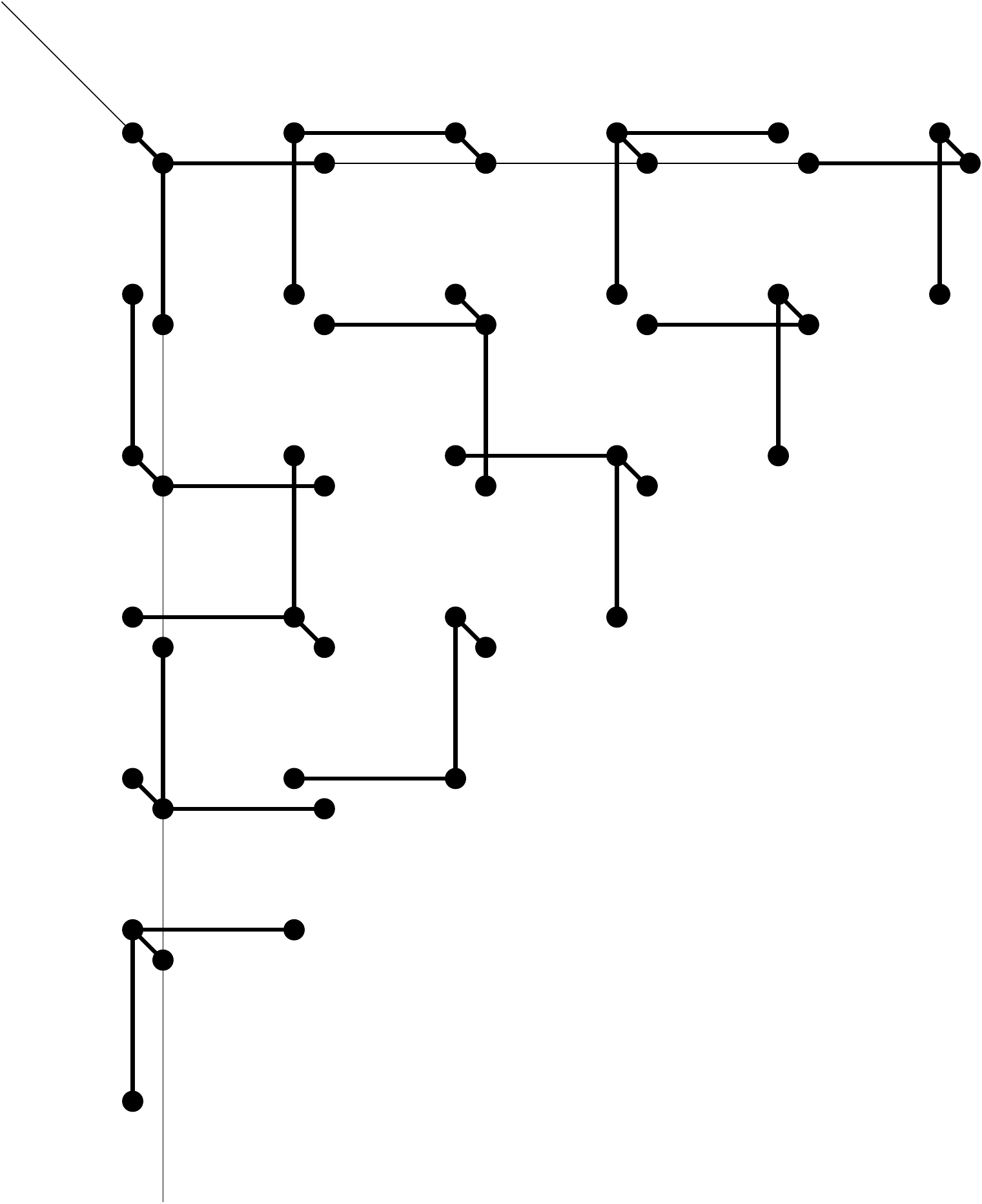}\label{fig:sl3_51}}
\subfigure[$(3+4k,1)$]{\includegraphics[scale=.3]{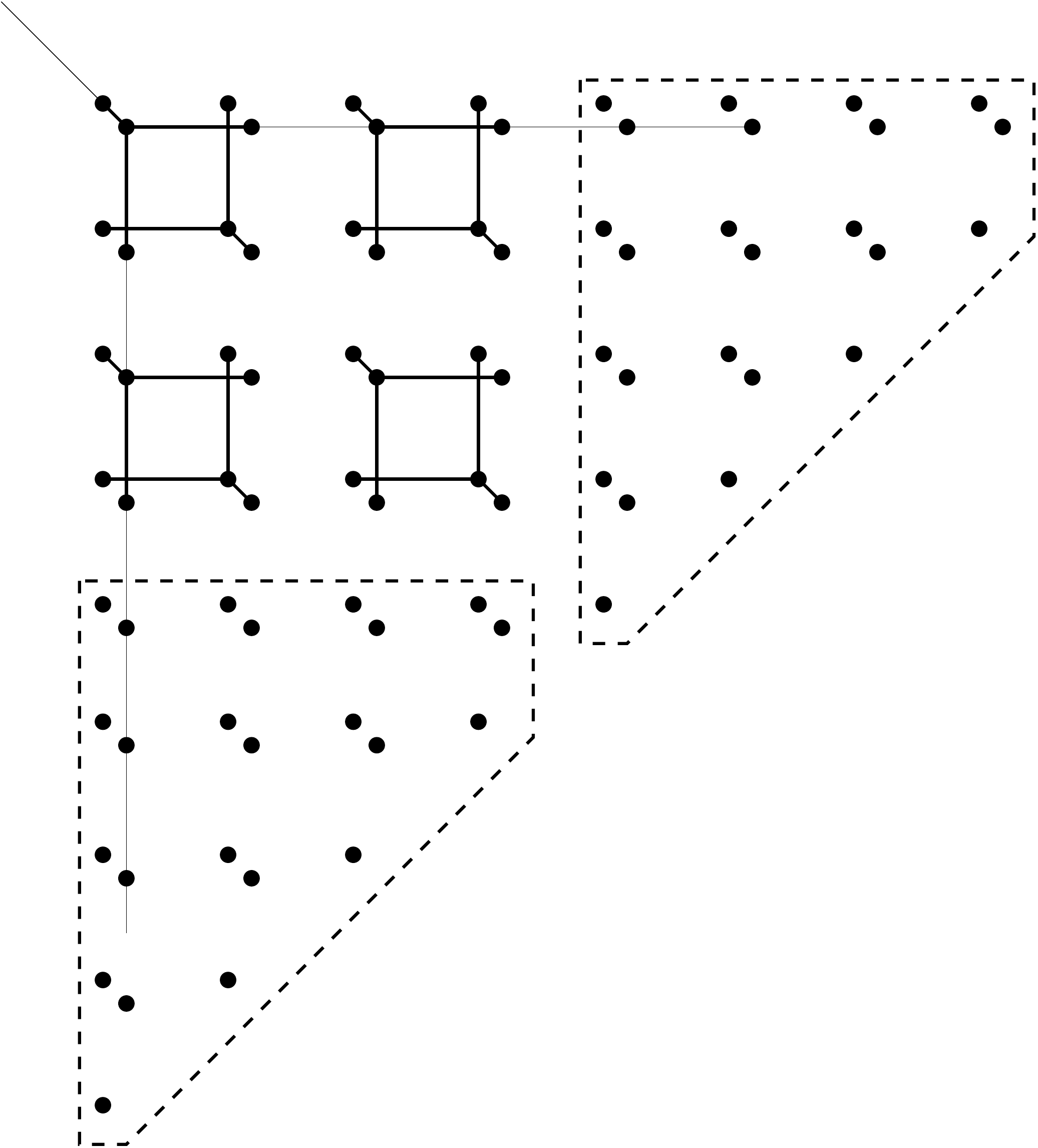}\label{fig:sl3_1ind2}}
\caption{Pictures for $(\text{odd},1)$}\label{fig:sl3_x1odd}
\end{center}
\end{figure}

In what follows we will assume that $m \geq n$ when convenient. We
first prove, by induction over $m$ non-defectiveness for $(m,1)$ and
$(m,2)$, and then do induction over $n$ to conclude the proof.  Figure
\ref{fig:sl3_11} for $(m,n)=(1,1)$ is {\em not} non-defective, reflecting
that the adjoint minimal orbit---the cone over which is the cone of $3
\times 3$-matrices with trace $0$ and rank $\leq 1$---is defective. Figure
\ref{fig:sl3_21}, however, shows a non-defective picture for $(2,1)$,
and from this picture one can construct non-defective pictures for
$(2+2k,1)$ by putting it to the right of $k$ pictures, each of which
consists of cubes and a single tetrahedron; Figure \ref{fig:sl3_1ind1}
illustrates this for the step from $(2,1)$ to $(2+2,1)$.

Figure \ref{fig:sl3_31} shows a non-defective picture for $(3,1)$, and
Figure \ref{fig:sl3_51} a non-defective picture for $(5,1)$. From these
we can construct non-defective pictures for $(3+4k,1)$ and $(5+4k,1)$,
respectively, by putting them to the right of $k$ pictures, each of which
consists of a few cubes plus a non-defective picture for $(3,1)$---Figure
\ref{fig:sl3_1ind2} illustrates this for the step from $(3,1)$ to $(7,1)$.
This settles $(n,1)$.

\begin{figure}
\begin{center}
\subfigure[$(2,2)$; defective]{\includegraphics[scale=.25]{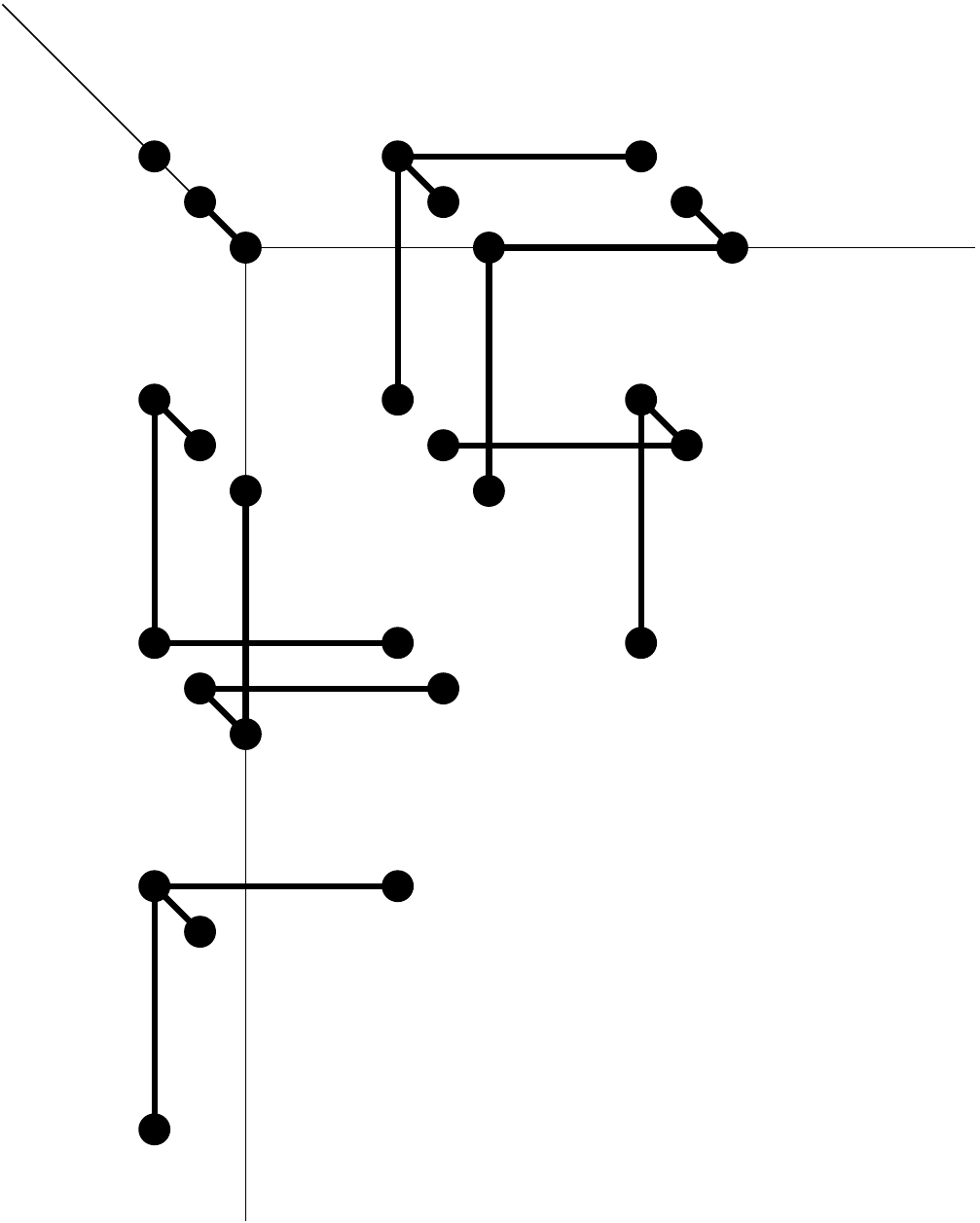}\label{fig:sl3_22}}
\subfigure[$(3,2)$]{\includegraphics[scale=.25]{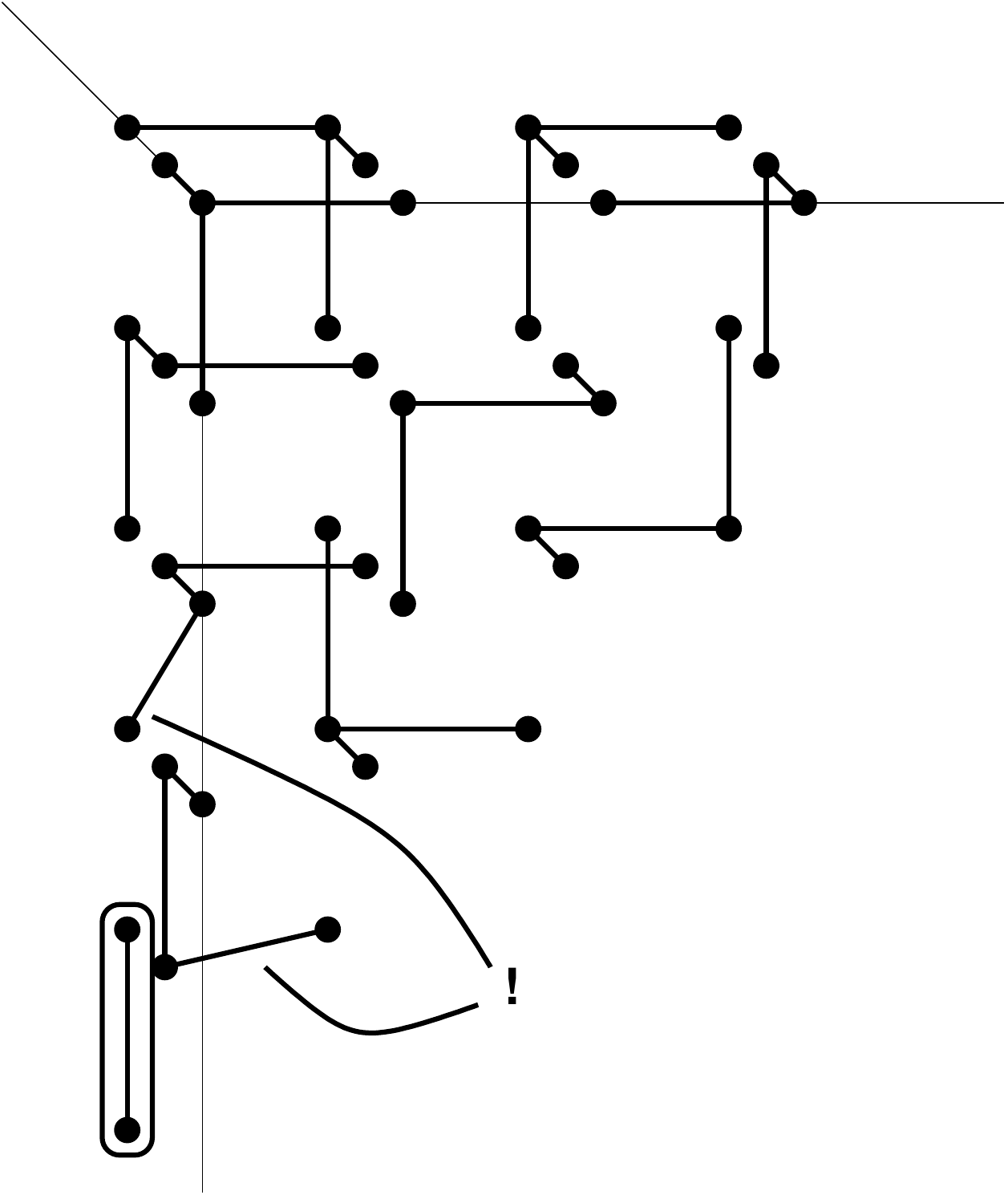}\label{fig:sl3_32}}
\subfigure[$(4,2)$]{\includegraphics[scale=.25]{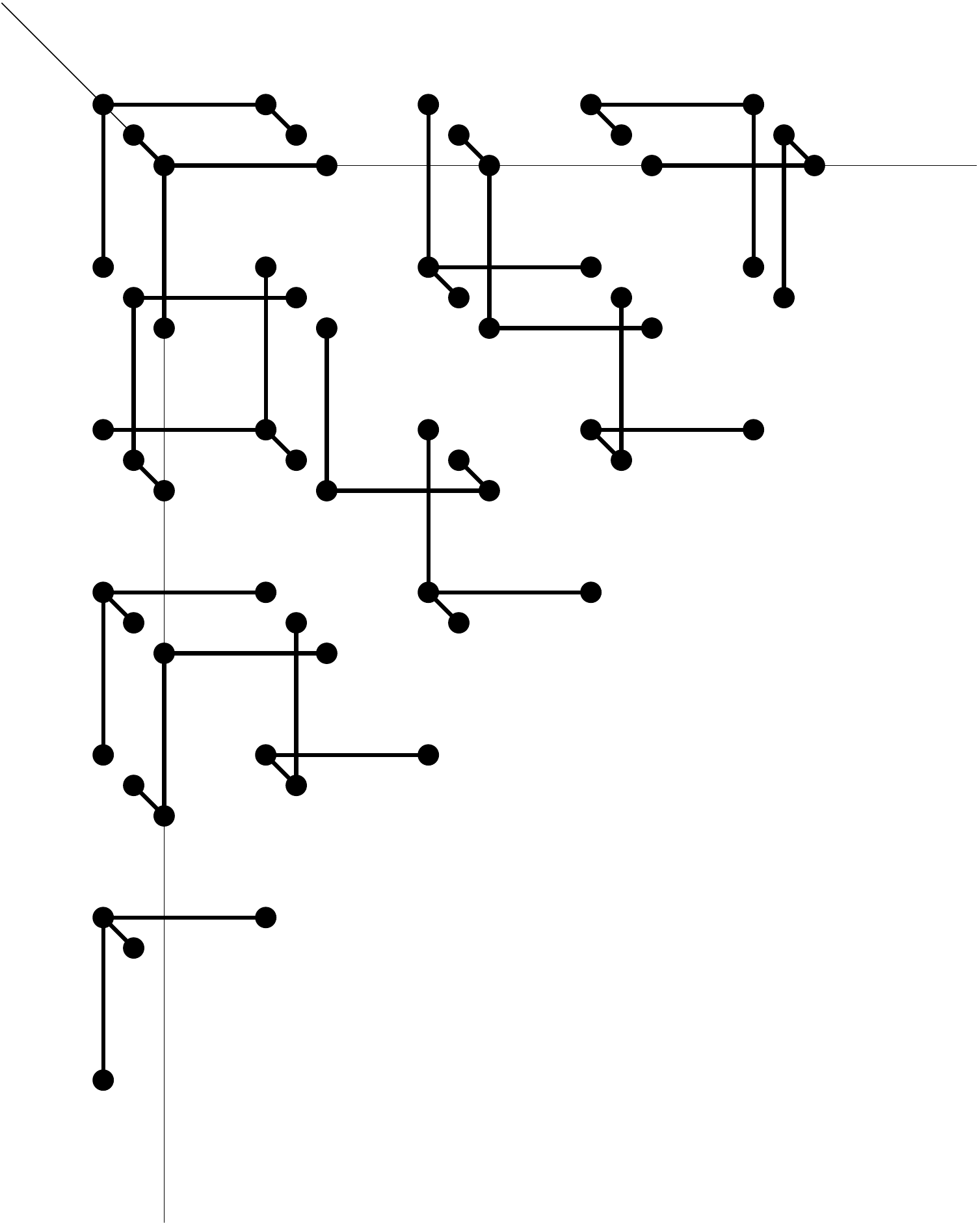}\label{fig:sl3_42}}
\subfigure[$(5,2)$]{\includegraphics[scale=.25]{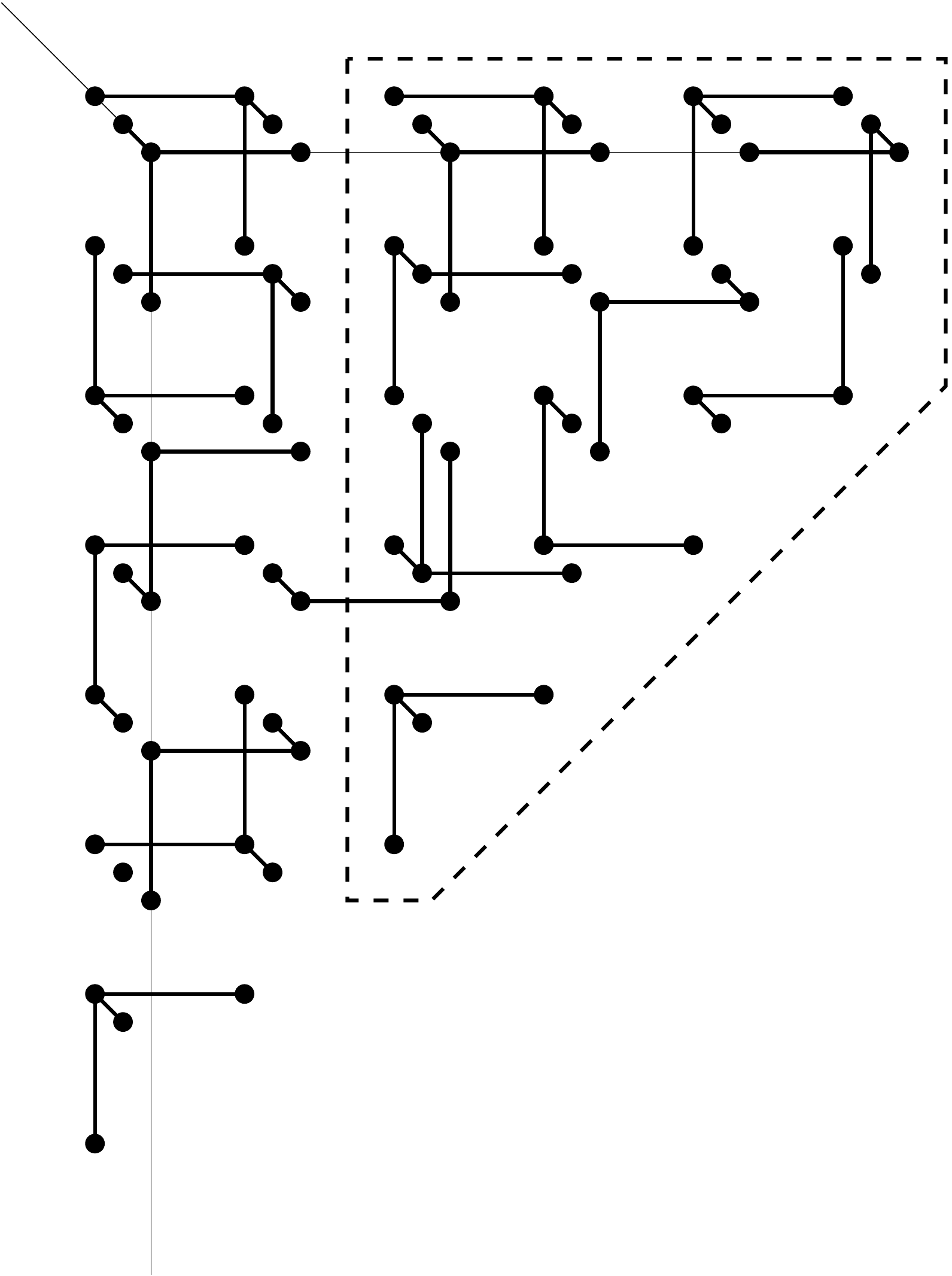}\label{fig:sl3_52}}
\subfigure[$(6,2)$]{\includegraphics[scale=.25]{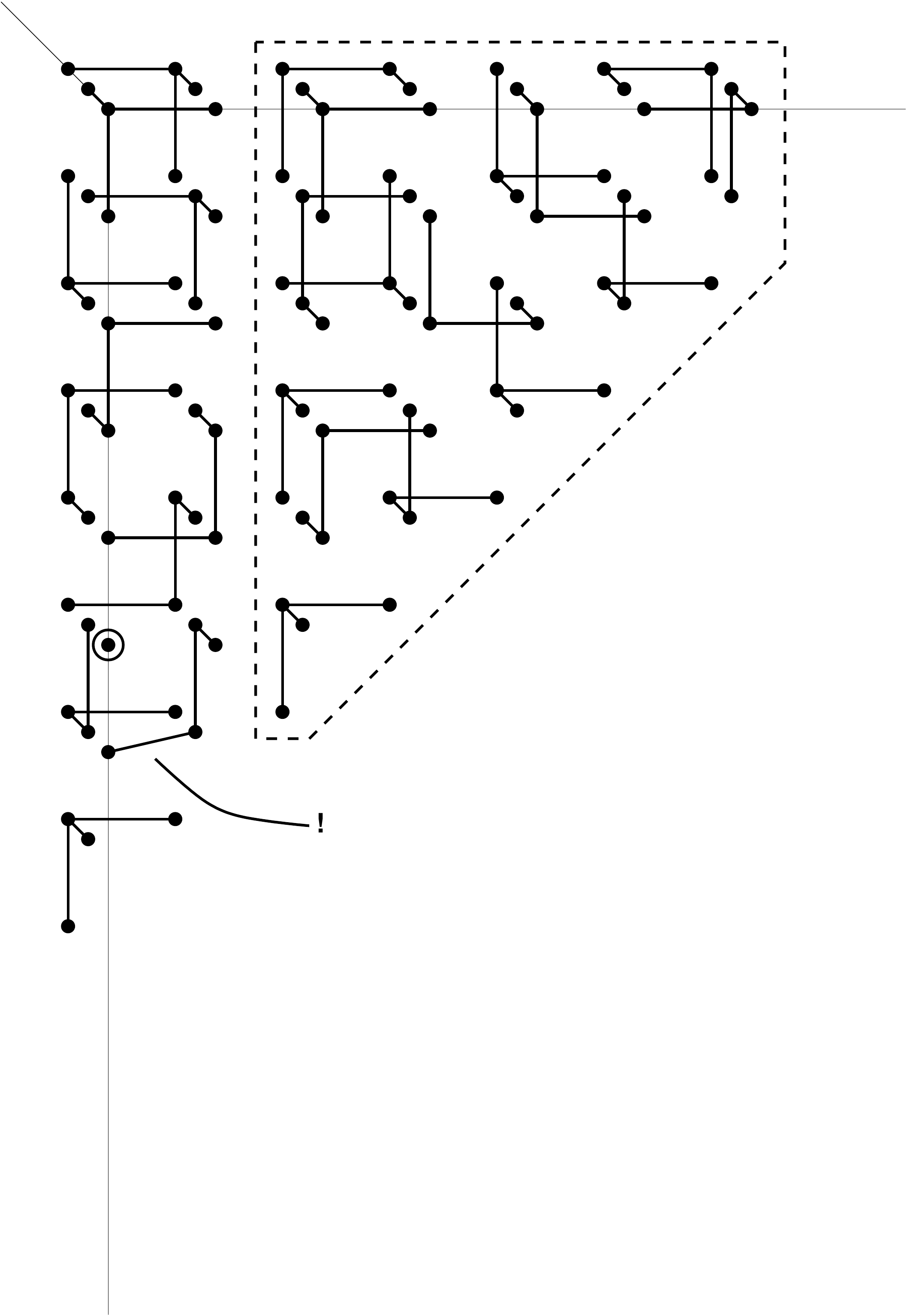}\label{fig:sl3_62}}
\subfigure[$(7,2)$]{\includegraphics[scale=.25]{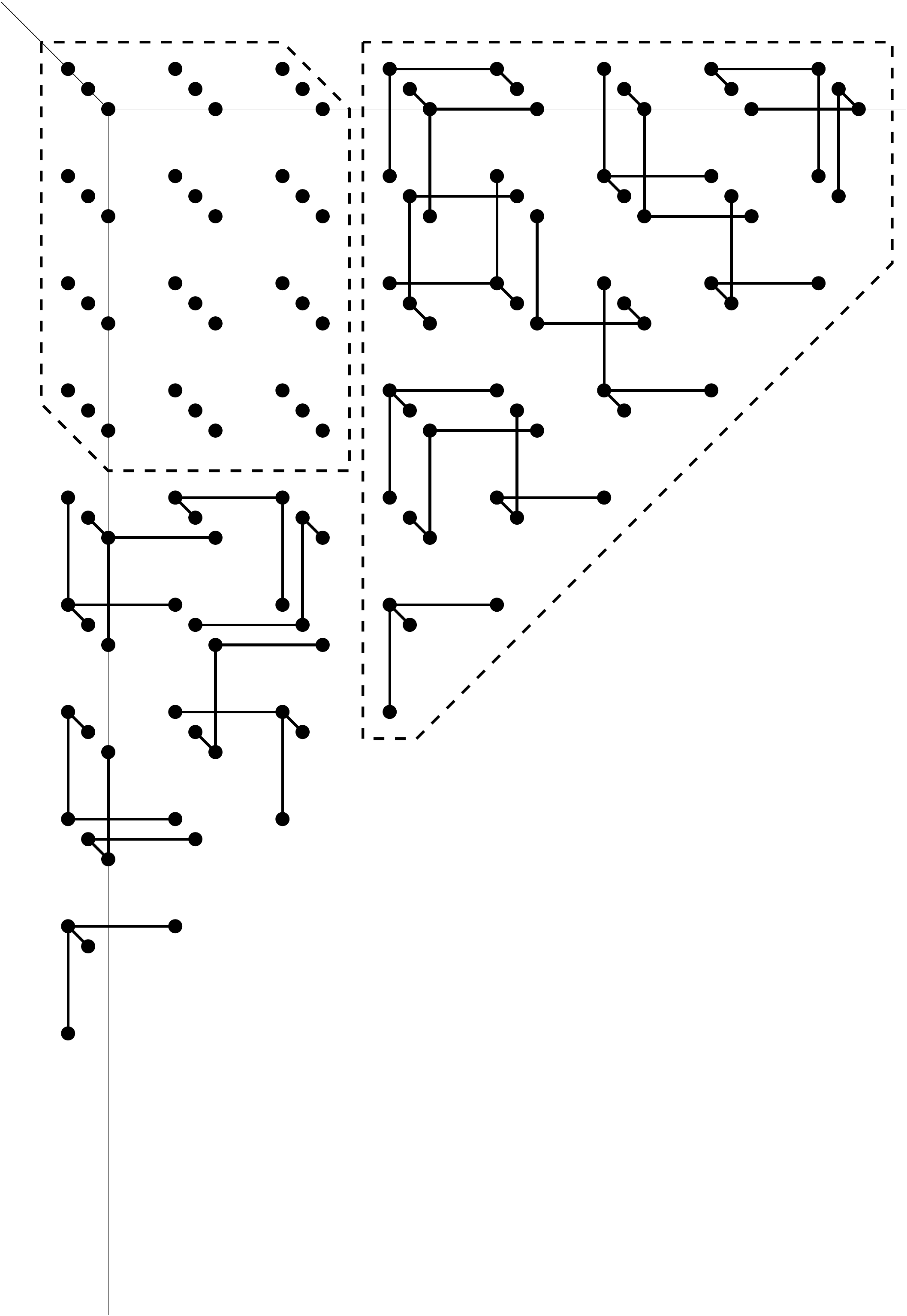}\label{fig:sl3_72}}
\subfigure[$(m+8,2)$]{\includegraphics[scale=.25]{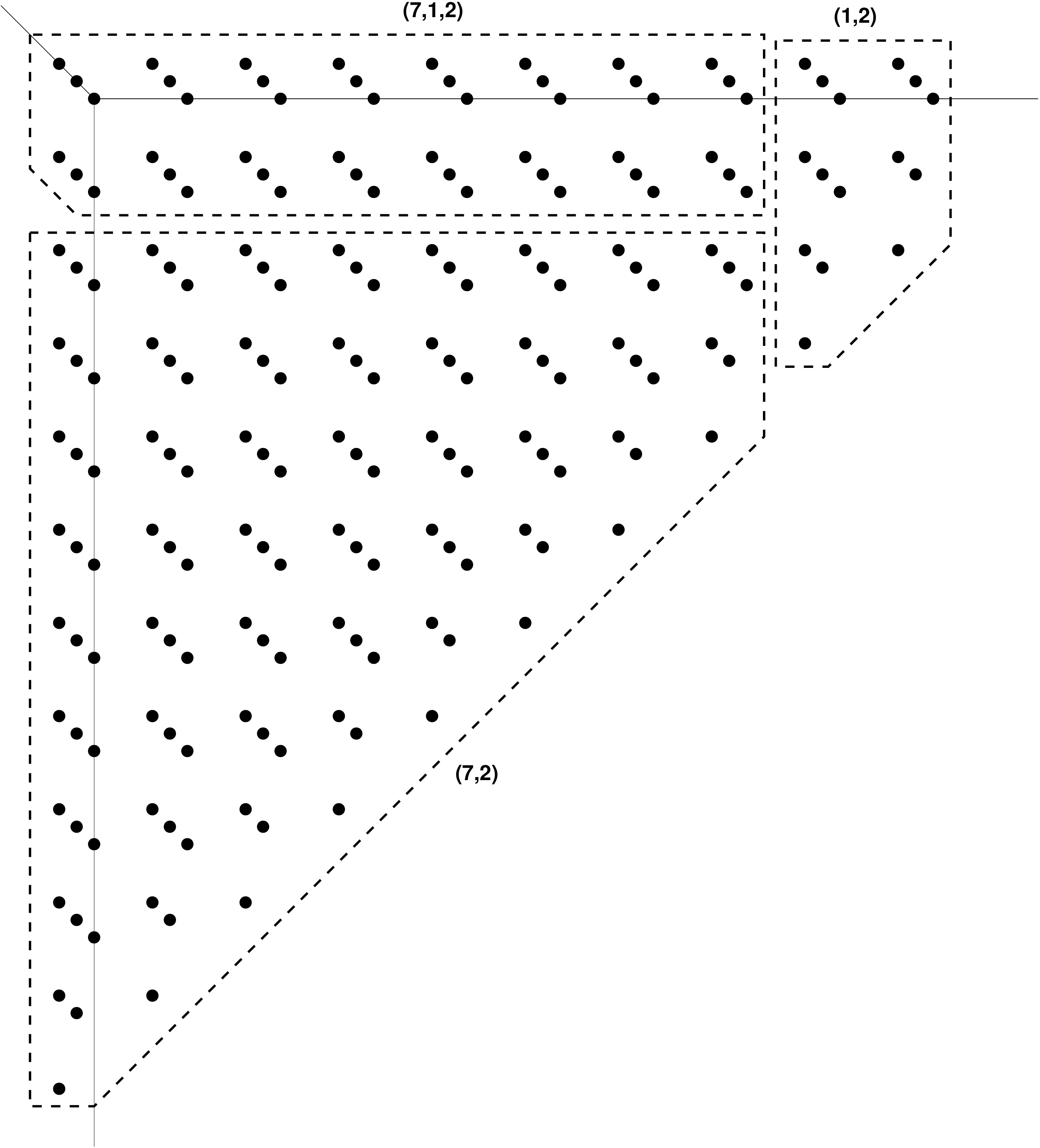}\label{fig:sl3_2ind}}
\caption{Pictures for $(*,2)$}\label{fig:sl3_x2}
\end{center}
\end{figure}

Figure \ref{fig:sl3_22} is defective: it reflects the fact that the $7$-th
secant variety of $X$ in the $(2,2)$-embedding has defect $1$.  Figure
\ref{fig:sl3_32} gives a non-defective picture for $(3,2)$. Note that
two non-standard cells $B_i$ are used here; this is because we will need
this picture in the picture for $(6,6)$.  Figure \ref{fig:sl3_42} gives
a non-defective picture for $(4,2)$. One can construct a non-defective
picture for $(5,2)$ based on the former, and a non-defective picture for
$(7,2)$ based on that for $(4,2)$; see Figure \ref{fig:sl3_52} and Figure
\ref{fig:sl3_72}. In the latter picture, one should fill in one copy of
our earlier non-defective picture for $\PP^1 \times \PP^1 \times \PP^1$
in its $(3,2,2)$-embedding. Figure \ref{fig:sl3_62} gives a non-defective
picture for $(6,2)$, based on that for $(4,2)$. A single cell $B_i$ is
non-standard, again for later use in the picture for $(6,6)$.  Similarly,
one can construct pictures for $(8,2)$ and $(10,2)$---which are left
out here because they take too much space.

Finally, from a non-defective picture for $(m,2)$ (with $m=1$ or $m
\geq 3$) one can construct a non-defective picture for $(m+8,2)$ by
inserting a picture consisting of a non-defective picture for $(7,2)$
and our non-defective picture for a $(7,m,2)$-block from the discussion
of $\PP^1 \times \PP^1 \times \PP^1$ in front---Figure \ref{fig:sl3_2ind}
illustrates this for $m=1$. This settles the cases where $m\geq n=2$.

\begin{figure}
\begin{center}
\includegraphics[scale=.5]{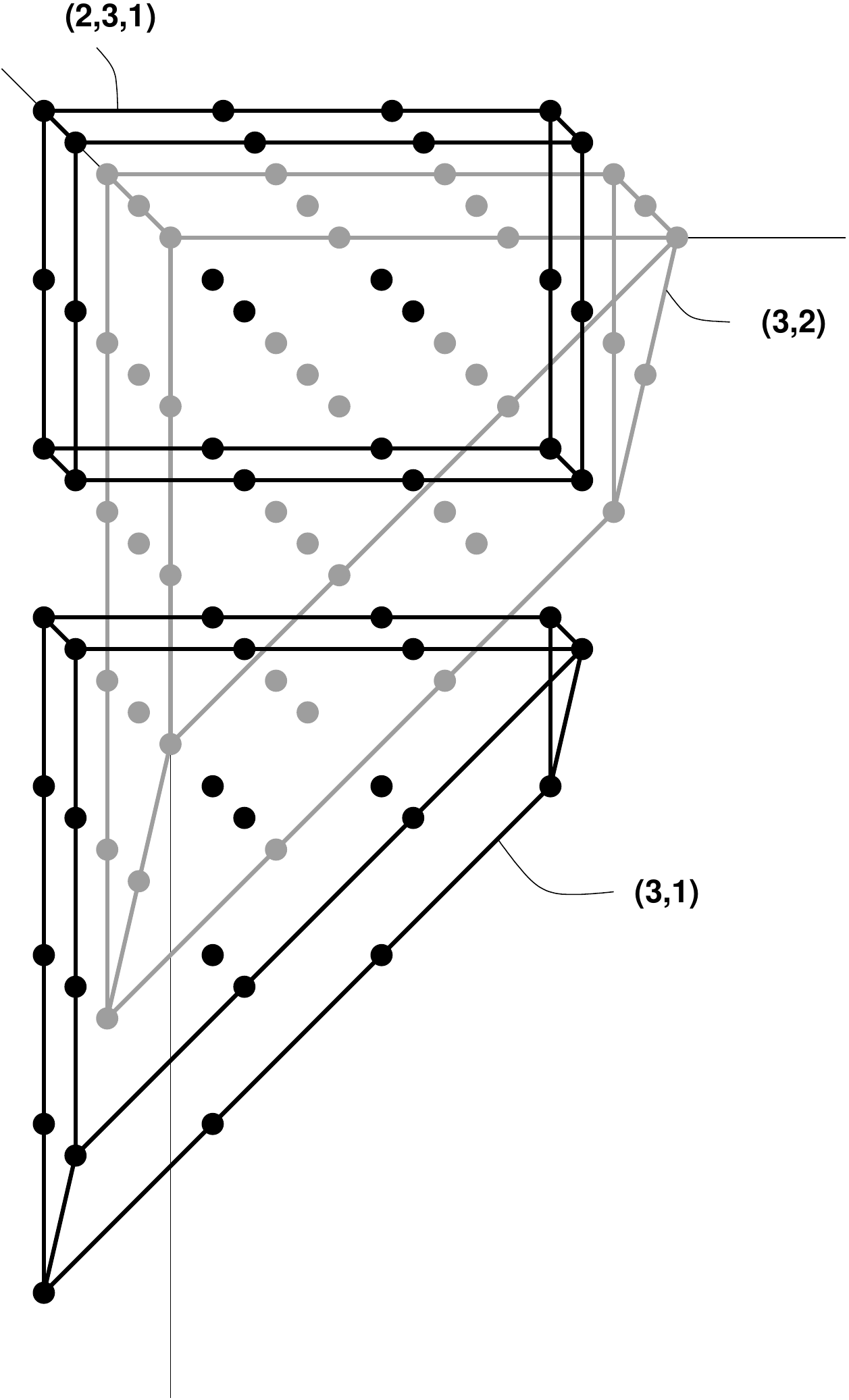}
\caption{Construction for $m$ or $n$ odd.}
\label{fig:sl3_oddind}
\end{center}
\end{figure}

Now all cases where at least one of $m$ and $n$ is odd can also be
settled. Indeed, suppose that $m,n \geq 3$ and that $m$ is odd. Write
$n+1=2q + (r+1)$ with $r \in \{1,2\}$. Then we can construct a
non-defective picture for $(m,n)$ by taking our non-defective picture
for $(m,r)$ and succesively stacking $q$ non-defective pictures of
two layers on top, each of which pictures with a number of vertices
divisible by $4$. These layers can be constructed as follows: the
$i$-th layer consists of our non-defective picture for $\PP^1 \times
\PP^1 \times \PP^1$ with parameters $(r+2i-2,m,1)$ (lying against the
$(n_2,n_3)$-plane) and a non-defective picture for $\mF$ with parameters
$(m,1)$. As $m$ is odd, each of these two blocks has a number of vertices
divisible by $4$. This construction is illustrated for $m=3$ and $n=4$
in Figure \ref{fig:sl3_oddind}, where one extra layer is put on top of
the ``ground layer''.

\begin{figure}
\begin{center}
\subfigure[$(4,4)$]{\includegraphics[scale=.25]{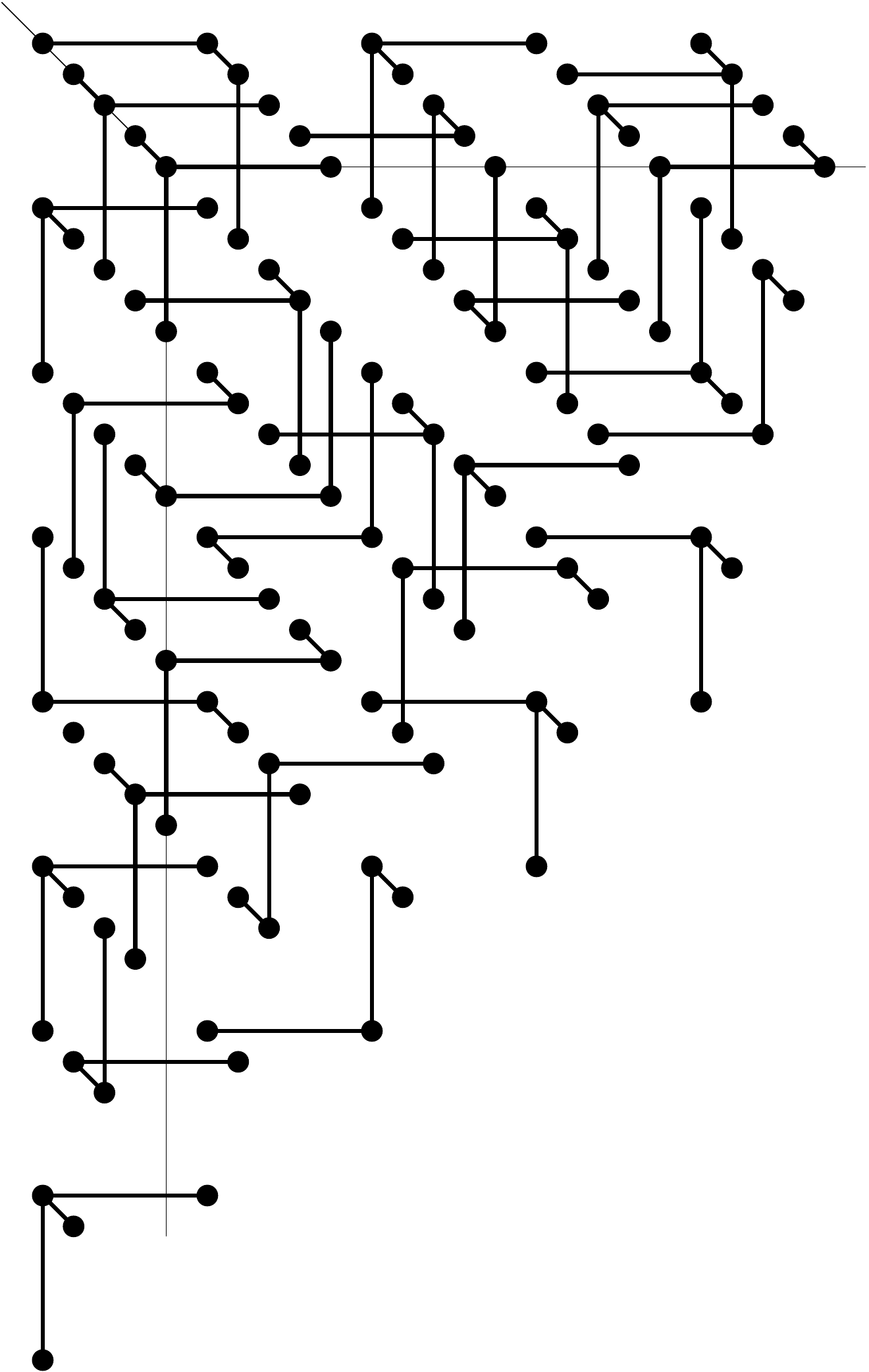}\label{fig:sl3_44}}
\subfigure[$(4,6)$]{\includegraphics[scale=.25]{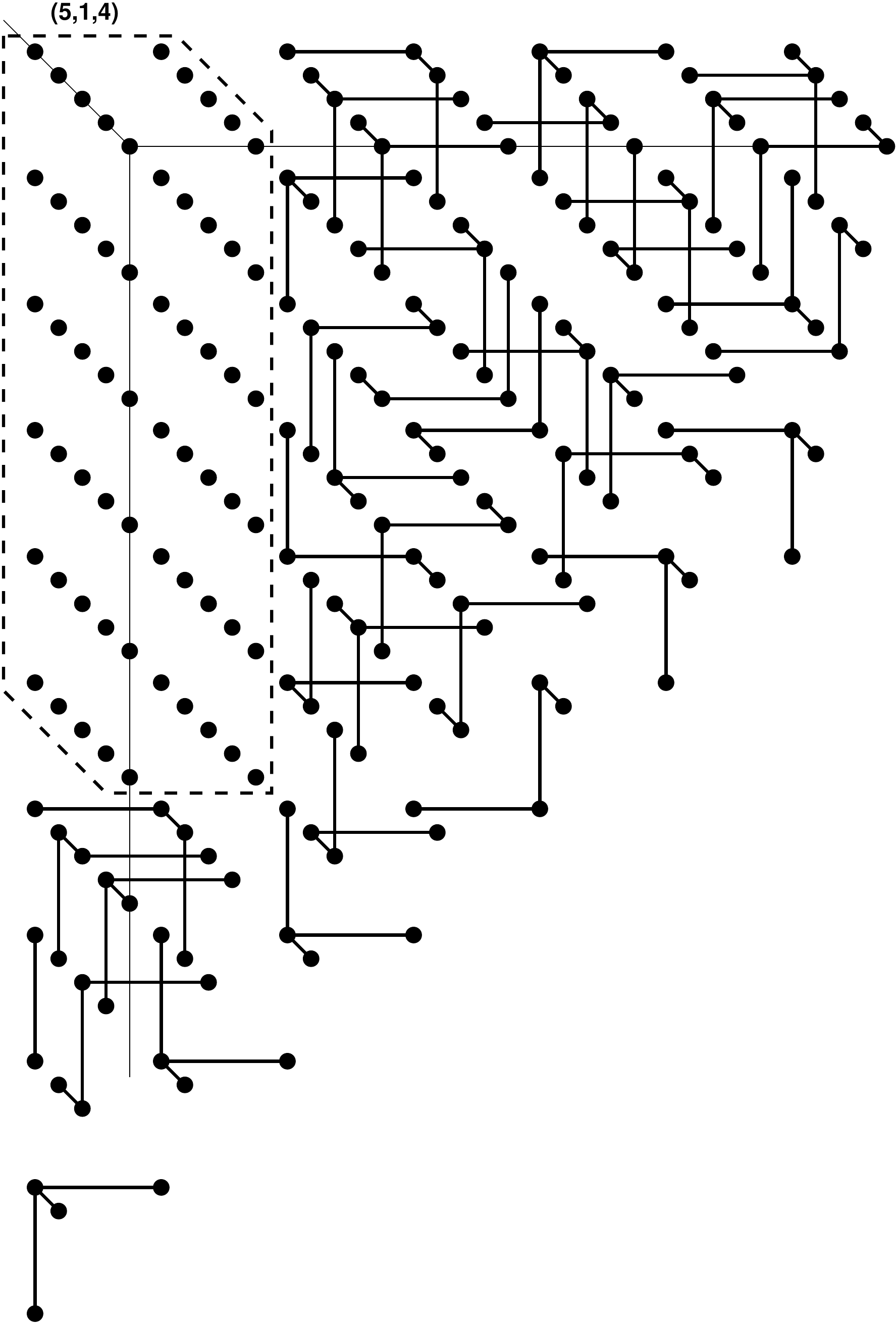}\label{fig:sl3_64}}
\subfigure[$(6,6)$]{\includegraphics[scale=.25]{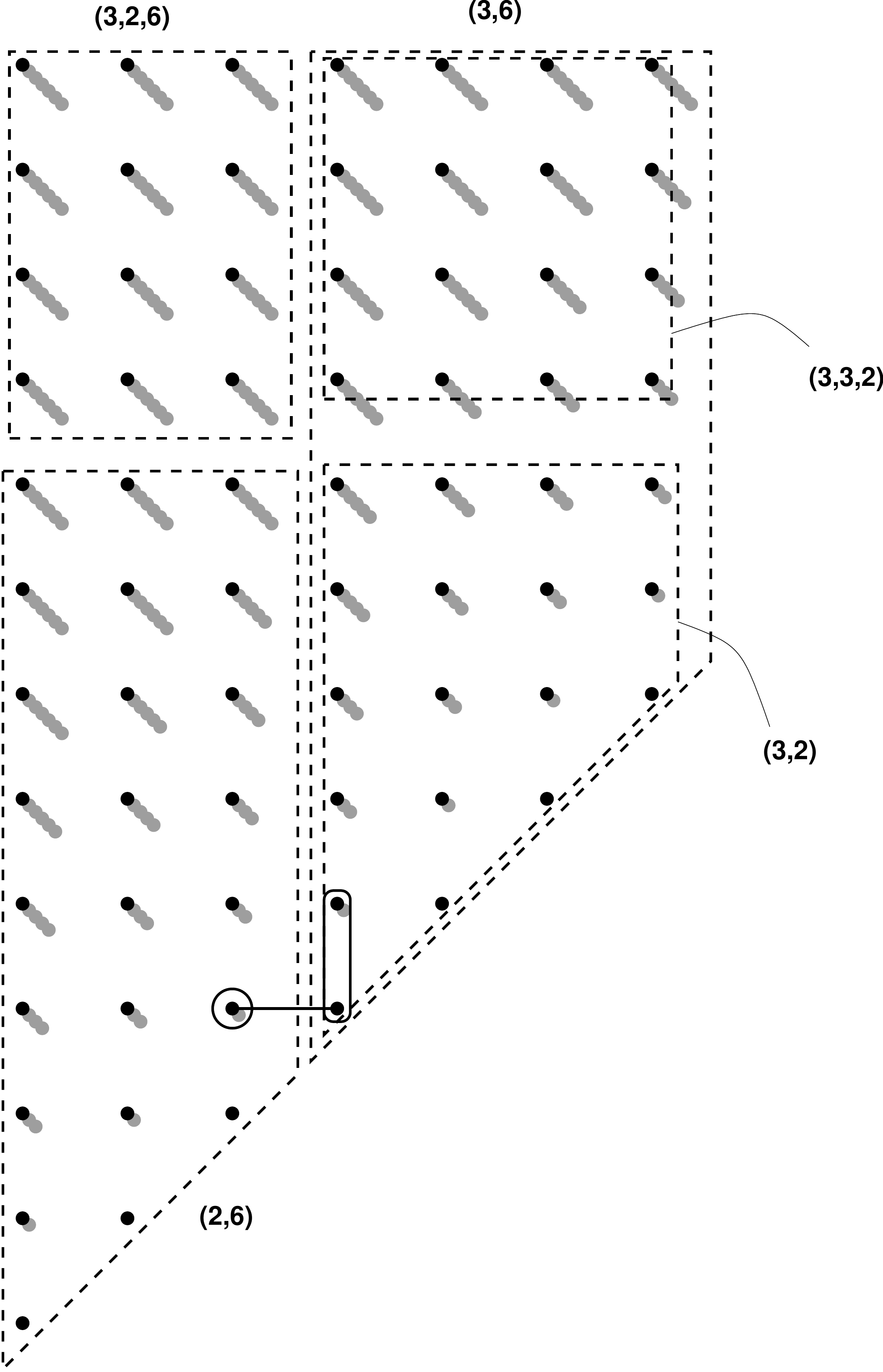}\label{fig:sl3_66}}
\end{center}
\end{figure}

Only the cases remain where $m$ and $n$ are both even. We first argue that
we can now reduce the discussion to a finite problem: if $m \geq 7$ and
$n \geq 2$, then we can compose a non-defective picture for $(m,n)$ from
one non-defective picture for $(7,n)$ (which exists by the above), one
non-defective $\PP^1 \times \PP^1 \times \PP^1$-picture for $(7,m-8,n)$
(both of these have numbers of vertices divisible by $4$), and one
non-defective picture for $(m-8,n)$---if such a picture exists. Hence we
may assume that $m<7$. Similarly, by using Remark \ref{re:Transposition}
we may assume that $n<7$. Using that $m,n$ are even, and that $m,n>2$
(which we have already dealt with), we find that only $(4,4)$, $(4,6)$
(or $(6,4)$), and $(6,6)$ need to be settled---as done in Figures
\ref{fig:sl3_44}--\ref{fig:sl3_66}. The picture for $(4,6)$ uses our
picture for the Segre-Veronese embedding of $\PP^1 \times \PP^1 \times
\PP^1$ of degree $(5,1,4)$. The picture for $(6,6)$ is built from a
picture for $\PP^1 \times \PP^1 \times \PP^1$ of weight $(3,2,6)$, and
pictures for $\mF$ of weights $(2,6)$ and $(3,6)$; the latter picture,
in turn, can be constructed as outlined above, except that, in order to
line up the single edge in $(3,6)$ and the single vertex in $(2,6)$, the
order of the building blocks for $(3,6)$ is altered: the $\mF$-picture
for $(3,2)$ comes on top, next to a $(\PP^1)^3$ picture for $(3,3,1)$,
and under these an $\mF$-picture for $(3,3)$. This concludes the proof
of Theorem \ref{thm:F}.

\end{document}